\documentclass[reqno,11pt]{amsart}    

\usepackage{latexsym,amsfonts,amsmath,amsthm,amssymb,euscript}
\usepackage[english]{babel}
\usepackage[latin1]{inputenc}
\usepackage[T1]{fontenc}

\newcounter{compteur}
\newtheorem{theop}[compteur]{Theorem}
\newtheorem{defi}{Definition}[section]

\newtheorem{lemm}[defi]{Lemma}
\newtheorem{prop}[defi]{Proposition}

\newtheorem{corr}[defi]{Corollary}

\newtheorem{hypo}{Assumption}

\DeclareMathSymbol{\leqslant}{\mathalpha}{AMSa}{"36}      
\DeclareMathSymbol{\geqslant}{\mathalpha}{AMSa}{"3E}      
\DeclareMathSymbol{\eset}{\mathalpha}{AMSb}{"3F}          
\renewcommand{\leq}{\;\leqslant\;}                        
\renewcommand{\geq}{\;\geqslant\;}                        

\numberwithin{equation}{section}

\newcommand{\petit}[1]{\mbox{\tiny  $#1$}}
\newcommand{\N}{\ensuremath{\mathbb{N}}}                                  
\newcommand{\Z}{\ensuremath{\mathbb{Z}}}                                  
\newcommand{\R}{\ensuremath{\mathbb{R}}}                                  
\newcommand{\D}{\ensuremath{\mathbb{D}}}                                  
\newcommand{\Indic}{\ensuremath{\mathbf{1}}}                              
\newcommand{\Pare}[1]{\ensuremath{\left( #1 \right)}}                     
\newcommand{\Prob}[1]{\ensuremath{\mathbf{P}\left(#1\right)}}             
\renewcommand{\P}{\ensuremath{\mathbf{P}}}                                
\newcommand{\Pdeb}{\mathbb{P}}
\newcommand{\Omegadeb}{\Omega}
\newcommand{\Omegabon}{\mathbf{\Omega}}
\newcommand{\Esp}[1]{\ensuremath{\mathbf{E}\left(#1\right)}}              
\newcommand{\EE}{\ensuremath{\mathbf{E}}}                                 
\newcommand{\Cste}[1]{\ensuremath{\mathbf{C}_{#1}}}                       
\newcommand{\Event}[1]{\ensuremath{\mathcal{E}_{#1}}}                     
\newcommand{\cadlag}{c\`adl\`ag}                                          
\newcommand{\V}{\ensuremath{\mathbb{V}}}                                  
\renewcommand{\S}{\ensuremath{\mathbb{S}}}                                
\newcommand{\Mittag}{\ensuremath{E}}                                      
\newcommand{\M}{\ensuremath{\mathbf{M}}}                                  
\renewcommand{\H}{\ensuremath{\mathbf{H}}}                                
\newcommand{\T}{\ensuremath{\mathbf{T}}}                                  
\newcommand{\RenH}{\ensuremath{\mathbf{R}}}                               
\newcommand{\WartW}{\ensuremath{\mathbf{W}}}                              
\newcommand{\Achange}{\ensuremath{\mathbb{A}}}                            
\newcommand{\Tchange}{\ensuremath{\mathbb{T}}}                            
\newcommand{\Rbes}{\ensuremath{\mathcal{R}}}                              
\newcommand{\Zbes}{\ensuremath{\mathcal{Z}}}                              

\begin{document}
\title[Diffusion in a random environment]{Limiting behavior of a diffusion in an asymptotically stable environment}
\author{Arvind Singh}
\maketitle
\begin{center}
\it{Laboratoire de Probabilités et Modèles Al\'eatoires,\\
 Université Pierre et Marie Curie,\\ 175 rue du Chevaleret,
75013
Paris, France.\\
e-mail: arvind.singh@ens.fr}
\end{center}
\bigskip
\begin{abstract}
Let $\V$ be a two sided random walk and let $X$ denote a real valued
diffusion process with generator
$\frac{1}{2}e^{\V_{[x]}}\frac{d}{dx}\Pare{e^{-\V_{[x]}}\frac{d}{dx}}$.
This process is known to be the continuous equivalent of the one
dimensional random walk in random environment with potential $\V$.
Hu and Shi (1997) described the L\'evy classes of $X$ in the case
where $\V$ behaves approximately like a Brownian motion. In this
paper, based on some fine results on the fluctuations of random walks
and stable processes, we obtain an accurate image of the almost
sure limiting behavior of $X$ when $\V$ behaves asymptotically like a
stable process. These results also apply for the corresponding random
walk in random environment.
\end{abstract}

\bigskip
\noindent{\bf Key words. } Random environment, stable process,
iterated logarithm law.

\bigskip
\noindent{\bf MSC 2000. } 60K37, 60J60, 60G52, 60F15.
\vspace*{0.8cm}
\section{Introduction}
Let $\Pare{\V_x,x\in\R}$ be a \cadlag, real-valued locally bounded
stochastic process on some probability space $(\Omegadeb,\Pdeb)$ with
$\V_0=0$ a.s. Let also $(X_t)_{t\geq 0}$ be the coordinate process
on the space of continuous functions $C([0,\infty))$ equipped with
the topology of uniform convergence on compact set and the
associated $\sigma$-field. For each realization of $\V$, let $P_\V$
be a probability on $C([0,\infty))$ such that $X$ is a diffusion
process with $X_0 = 0$ and generator
\begin{equation*}
\frac{1}{2}e^{\V_{x}}\frac{d}{dx}\Pare{e^{-\V_{x}}\frac{d}{dx}}\hbox{.}
\end{equation*}
It is well known that $X$ may be constructed from a standard
Brownian motion through a change of scale and a change of time \cite{Ito1}.
We consider the annealed probability $\P$ on 
$\Omegabon =\Omegadeb\times C([0,\infty))$ defined as the semi-direct product
$\P = \Pdeb\times P_\V$. $X$ under $\P$ is called a diffusion in the
random potential $\V$. This process was first studied by
Schumacher \cite{Schum1} and Brox \cite{Brox1} who proved that when
$\V$ is a Brownian motion $X_t/\log^2 t$ converges in law, as $t$ goes
to infinity, to some non degenerate distribution on $\R$. Extension
of this result when $\V$ is a stable process may be found in
\cite{Schum1,Tanak1,Cheli1}. In this paper, we are concerned with
the case where $\V$ is a two sided random walk. More precisely,
$(\V_x\hbox{ , }x\in\R)$ satisfies:
\begin{equation*}
\left\{%
\begin{array}{l}
\V\hbox{ is identically $0$ on $(-1,1)$,} \\
\V\hbox{ is flat on $(n,n+1)$ for all $n\in\Z$,} \\
\V\hbox{ is right continuous on $[0,\infty)$ and left continuous on $(-\infty,0]$,} \\
(\V_{n+1}-\V_n)_{n\in\Z}\hbox{ is a sequence of i.i.d. variables
under \P}.
\end{array}%
\right.
\end{equation*}
Our goal is to describe the almost sure asymptotics of $X_t$,
$\sup_{s\leq t}X_s$ and $\sup_{s\leq t}|X_s|$. This has been done by Hu and Shi \cite{HuShi1}
in the case where $\V$ behaves roughly like a Brownian motion. We
will instead consider the more general setting where a typical step
of the random walk is in the domain of attraction of a stable law.
In fact, we will make an assumption similar to that of Kawazu,
Tamura and Tanaka \cite{Tanak1}, that is, in all the following:
\begin{hypo} \label{hyp1} There exists a positive  sequence
$\Pare{a_n}_{n\geq 0}$ such that
\begin{equation*}
\frac{\V_n}{a_n}\underset{n\to\infty}{\overset{\hbox{law}}{\longrightarrow}}\S
\end{equation*}
where $\S$ is a random variable whose law is strictly stable with index $\alpha \in
(0,2]$ and whose density is everywhere positive on $\R$.
\end{hypo}
This implies of course that $\V_{-n}/a_n$ converges in law toward
$-\S$. It is known that the norming sequence $(a_n)$ is regularly
varying with index $1/\alpha$ and we can without loss of
generality assume that $(a_n)$ is strictly increasing with $a_1 =
1$. We will denote by $a\Pare{\cdot}$ a continuous, strictly
increasing interpolation of $(a_n)$ and $a^{-1}\Pare{\cdot}$ will
stand for its inverse. It is to be noted that $a\Pare{\cdot}$ and
$a^{-1}\Pare{\cdot}$ are respectively regularly varying with index $1/\alpha$ and
$\alpha$. Let $p$  denote the positivity parameter of $\S$ and $q$
its negativity parameter, namely:
\begin{equation*}
 p = \Prob{\S > 0} = 1 - \Prob{\S< 0} = 1-q.
\end{equation*}
The assumption that $\S$ has a positive density in the whole of $\R$
implies that $p,q \in (0,1)$. More precisely for $\alpha
>1$ it is known  \cite{Zolot1} that $1-1/\alpha \leq p,q \leq
1/\alpha$. In any case:
\begin{equation*}
0 < \alpha p, \alpha q \leq 1.
\end{equation*}
Note also that the  Fourier transform of $\S$ is well known to be:
\begin{equation} \label{caractS}
\Esp{e^{i\lambda \S}} =
e^{-\gamma|\lambda|^{\alpha}\Pare{1-i\frac{\lambda}{|\lambda|}\tan\Pare{\pi\alpha\Pare{p-\frac{1}{2}}}}}
\end{equation}
where $\gamma$ is some strictly positive constant. Let us now extend
$\S$ into a two sided strictly stable process $(\S_x\hbox{ , }x\in\R)$
such that $\S_1$ has same law as $\S$. By two sided, we mean that
the processes $(\S_t\hbox{ , }t\geq 0)$ and $(-\S_{-t}\hbox{ ,
}t\geq 0)$ are independent, both \cadlag, and have the same law.
Notice in particular that, when $\alpha = 1$, $\S$ is a symmetric
Cauchy process with drift, whereas for $\alpha=2$ we have $p=1/2$
and $\S$ is a Brownian motion. Furthermore,  the extremal cases
$\alpha p =1$ (resp. $\alpha q = 1$) can only happen when $\alpha
>1$ and are equivalent to the assumption that $\S$ has no positive
jumps (resp. no negative jumps). When $\S$ has no positive jumps, it
is known that the Fourier transform can be extended such that
\begin{equation} \label{ee1}
\EE\big(e^{\lambda \S_1}\big) = e^{\gamma' \lambda^\alpha} \hbox{for all
$\lambda \geq 0$}
\end{equation}
where $\gamma'$ is a positive constant that we will assume to be $1$
(we can reduce to this case by changing the norming sequence $a_n$).
Similarly, when $\S$ has no negative jumps, we will assume
$\Esp{\exp(-\lambda \S_1)} = \exp(\lambda^\alpha)$ for all $\lambda
\geq 0$. Let $\Mittag_\alpha$ denote the Mittag-Leffler function
with parameter $\alpha$:
\begin{equation*}
\Mittag_\alpha(x) = \sum_{n=0}^{\infty}\frac{x^{n}}{\Gamma(\alpha n + 1)} \hbox{ for $x\in\R$}.
\end{equation*}
Let also define $-\rho_1(\alpha)$ to be the first negative root of
$\Mittag_\alpha$ and $-\rho_2(\alpha)$ to be the first negative root
of $\alpha x \Mittag_\alpha''(x) + (\alpha-1)\Mittag_\alpha'(x)$.
The first result of this paper is a law of the iterated logarithm
for the limsup of the diffusion in random environment $X$.
\begin{theop} \label{MainTheo1}Under the annealed probability
$\P$, almost surely:
\begin{equation*}
\limsup_{t\to\infty}\frac{X_t}{a^{-1}\Pare{\log t}\log\log\log t} =
\frac{1}{K^\#}
\end{equation*}
where $K^\# \in (0,\infty)$ is a constant that only depends on
the limit law $\S$ and is given by the formula:
\begin{equation*}
K^\# = -\lim_{t\to\infty}\frac{1}{t}\log \P\Big(\sup_{0\leq u\leq
v\leq t}\Pare{\S_v - \S_u} \leq 1\Big).
\end{equation*}
Furthermore, when $\S$ is completely asymmetric, the value of $K^\#$ is given by:
\begin{equation*}
K^\# =
\left\{%
\begin{array}{l}
    \rho_1 (\alpha) \hbox{ when $\S$ has no positive jumps}, \\
    \rho_2 (\alpha) \hbox{ when $\S$ has no negative jumps}. \\
\end{array}%
\right.
\end{equation*}
\end{theop}
Note that $X_t$ and $\sup_{s\leq t}X_s$ have the same upper functions, hence
Theorem \ref{MainTheo1} also holds with $\sup_{s\leq t}X_s$ in place of $X_t$.
From a symmetry argument:
\begin{equation*}
\limsup_{t\to\infty}\frac{-\inf_{s\leq t}X_t}{a^{-1}\Pare{\log
t}\log\log\log t} = \frac{1}{\widetilde{K}^\#} \hbox{ a.s.}
\end{equation*}
where $\widetilde{K}^\# = -\lim_{t\to\infty}\log \Prob{\sup_{0\leq
u\leq v\leq t}\Pare{\S_{-v} - \S_{-u}} \leq 1}/t$, hence
\begin{equation*}
\limsup_{t\to\infty}\frac{\sup_{s\leq t}|X_t|}{a^{-1}\Pare{\log
t}\log\log\log t} = \frac{1}{\widetilde{K}^\# \wedge K^\#} \hbox{
a.s.}
\end{equation*}
In the case where $\alpha=2$, we have $\Mittag_\alpha (-x) = \cos(\sqrt{x})$ for all $x\geq 0$,
therefore $\widetilde{K}^\# = K^\# = \pi^2/4$ and we recover the law
of the iterated logarithm of Theorem $1.6$ of \cite{HuShi1}.

Let us denote $\T_n$ the $n^{\hbox{\tiny{th}}}$ strict descending
ladder index of the random walk $\V$, formally:
\begin{equation*}
\left\{%
\begin{array}{l}
    \T_0 =0, \\
    \T_{n+1} = \min\Pare{k>\T_n \hbox{ , } \V_{k} < \V_{\T_n}}.
\end{array}%
\right.
\end{equation*}
Since $\V$ is oscillatory, $\T_n$ is proper for all $n$.
Theorem $4$ of Rogozin \cite{Rogoz1} states that $\T_1$ is in the domain of
attraction of a positive stable law with index $q$, moreover $\T_1$ is in the
domain of normal attraction of this distribution if and only if
\begin{equation}\label{ro1}
\sum_{n=1}^{\infty}\frac{\Prob{\V_n<0} - q}{n} < \infty.
\end{equation}
Let $(b_n)$ denote a (strictly increasing) sequence of norming
constants for $\T_1$ and $b\Pare{\cdot}$ will stand for a continuous,
strictly increasing interpolation of this sequence. The function
$b^{-1}\Pare{\cdot}$ is therefore regularly varying with index $q$. The next
theorem characterizes the liminf behavior of $\sup_{s\leq t}X_s$.
\begin{theop} \label{MainTheo2} For any positive, non decreasing function $f$ we have:
\begin{equation*}
\P\Pare{\sup_{s\leq t} X_s \leq f(t) \hbox{ i.o. }} =
\left\{%
\begin{array}{l}
    0 \\
    1
\end{array}%
\right. \Longleftrightarrow
\int^{\infty}\frac{b^{-1}\Pare{f(t)}dt}{b^{-1}\Pare{a^{-1}\Pare{\log
t}}t\log t}
\left\{%
\begin{array}{l}
    <\infty \\
    =\infty.
\end{array}%
\right.
\end{equation*}
In particular, with probability $1$:
\begin{equation*}
\liminf_{t\to\infty}\frac{\Pare{\log\log
t}^{\beta}}{a^{-1}\Pare{\log t}}\sup_{s\leq t}X_s =
\left\{%
\begin{array}{ll}
    0, & \hbox{if $\beta<1/q$}, \\
    \infty, & \hbox{if $\beta>1/q$}. \\
\end{array}%
\right.
\end{equation*}
\end{theop}
Note that $(\ref{ro1})$ hold whenever $\V_1$ is strictly stable or
when $\Esp{\V_1^2} < \infty$ (according to Theorem $1$ of
\cite{Felle1}, p 575). In those cases, $\V_1$ is also in the domain
of normal attraction of $\S$ so that we can both choose
$a(x)=x^{1/\alpha}$ and $b(x) = x^{1/q}$ and the last theorem is
simplified:
\begin{equation*}
\P\Pare{\sup X_s \leq \frac{\log^\alpha t}{f(t)} \hbox{ i.o. }} =
\left\{%
\begin{array}{l}
    0 \\
    1
\end{array}%
\right. \Longleftrightarrow
\int^{\infty}\frac{dt}{tf^q(t)\log t}
\left\{%
\begin{array}{l}
    <\infty \\
    =\infty.
\end{array}%
\right.
\end{equation*}
In particular, the critical case $\beta = 1/q$ gives
\begin{equation*}
\liminf_{t\to\infty}\frac{(\log\log t)^{1/q}}{\log^\alpha
t}\sup_{s\leq t}X_s = \infty \hbox{ a.s.}
\end{equation*}

We are also interested in the asymptotic behavior of the bilateral
supremum $\sup_{s\leq t}|X_s|$. We already mentioned that the limsup
behavior of this process may be deduced from Theorem
\ref{MainTheo1}. Although we were not able to deal with the general
case (as it seems that many different behaviors may occur in the
completely asymmetric case, depending on the distribution tail of
$\V_1$) we can still obtain, when the limiting process has jumps
of both signs, an iterated logarithm law:
\begin{theop}\label{MainTheo2b}
When the limiting stable process $\S$ has jumps of both signs, we have,
for any increasing positive function $f$:
\begin{equation*}
\P\Pare{\sup_{s\leq t} |X_s| \leq \frac{a^{-1}(\log t)}{f(t)} \hbox{ i.o. }} =
\left\{%
\begin{array}{l}
    0 \\
    1
\end{array}%
\right. \Longleftrightarrow
\int^{\infty}\frac{dt}{tf(t)^2\log t}
\left\{%
\begin{array}{l}
    <\infty \\
    =\infty.
\end{array}%
\right.
\end{equation*}
In particular, with probability $1$:
\begin{equation*}
\liminf_{t\to\infty}\frac{\Pare{\log\log
t}^{\beta}}{a^{-1}\Pare{\log t}}\sup_{s\leq t}|X_s| =
\left\{%
\begin{array}{ll}
    0, & \hbox{if $\beta\leq 1/2$}, \\
    \infty, & \hbox{if $\beta>1/2$}. \\
\end{array}%
\right.
\end{equation*}
\end{theop}
Note that in this case, the limiting behavior does not depend on the
symmetry parameter and note also that this behavior is quite
different from the Brownian case (Theorem $1.7$ of \cite{HuShi1}).
This may be informally explained from the facts that when the
limiting process has jumps of both signs, typical valleys for the
diffusion are much deeper than in the Brownian case.

Although we are mainly concerned with the almost-sure behavior of $X$,
our approach also allows us to prove a convergence in law for the supremum process.
\begin{theop} \label{MainTheo3}
There exists a non degenerate random variable $\Xi$
 depending only on the limiting process $\S$ such that under the
annealed probability $\P$:
\begin{eqnarray*}
\frac{\sup_{s\leq t}X_s}{a^{-1}\Pare{\log
t}}\underset{t\rightarrow\infty}{\overset{\hbox{law}}{\longrightarrow}}\Xi.
\end{eqnarray*}
Moreover, when $\S$ has no positive jumps the law of $\Xi$ is
characterized by its Laplace transform
\begin{equation*}
\Esp{e^{-q\Xi}} = \Gamma\Pare{\alpha + 1}\frac{\Mittag'_\alpha(q)}{
\Mittag_\alpha(q)},
\end{equation*}
and in the case where $\S$ has no negative jumps:
\begin{equation*}
\Esp{e^{-q\Xi}} = (\alpha-1)\frac{\Mittag'_\alpha(q)}{\alpha q
\Mittag''_\alpha(q) + (\alpha-1)\Mittag'_\alpha(q)}.
\end{equation*}
\end{theop}
This paper is organised as follows: in section $2$, we prove
sharp results on the fluctuations of the potential $\V$ as well as
on the limiting stable process $\S$. These estimates which may be
of independent interest ultimately play an important role in the proof 
of the main theorems. In section $3$, we reduce the study of the hitting times
of $(X_t)$ to the study of some functionnals of the potential process $\V$.
This step is similar to \cite{HuShi1}, namely, we make use of Laplace's method
and the reader may refer to \cite{Shi1} for an overview of the key ideas. 
The proof of the main theorems are given in section $4$. We shall eventually
discuss these results in the last section, in particular, we show that
Theorems $\ref{MainTheo1}-\ref{MainTheo3}$ still hold when $\V$ is a stricly
stable process. We also explain how similar results can be deduced for a 
random walk in a random environment with an asymptotically stable potential.

\section{Fluctuations of $\V$ and $\S$}
In this section we prove several results about
fluctuations of the random walk $\V$. Some of these estimates will be
obtained via the study of the limiting process $\S$. In the first
subsection, we recall elementary properties of the stable
process $\S$ as well as a result of functional convergence of the
random walk toward the limiting stable process. In the following,
for any process $Z$, we will use indifferently the notation $Z_x$ or
$Z(x)$.
\subsection{Preliminaries and functional convergence in \D}\label{subsectionD}

We introduce the space $\D(\R_+ , \R)$ of \cadlag \hbox{ }functions
$Z:\R_+\to\R$ equipped with the Skorohod topology. Let $\theta$ stand
for the shift operator that is for any $Z\in\D(\R_+,\R)$ and any $x_0\geq 0$:
\begin{equation}
\Pare{(\theta_{x_0} Z)_x , x\geq 0} = \Pare{Z_{x+x_0}-Z_{x_0} , x\geq 0}
\end{equation}
Since our processes are double-sided, we will also need the space $\D(\R,\R)$
of functions $f:\R\to\R$ which are right continuous with left limits
on $[0,\infty)$ and left continuous with right limits on
$(-\infty,0]$ considered jointly with the associated Skorohod
topology. Recall that $\S$ and $\V$ have paths on $\D(\R,\R)$. We
will be interested in the following functionals: for any $a\in\R$ and
for any $Z\in\D(\R,\R)$ we define (we give two notations for each
definition):
$$
\begin{array}{rcccl}
\overline{Z}_a &=& F_a^{(1)}(Z) &=&
\left\{%
\begin{array}{ll}
    \sup_{y\in[0,a]} Z_y, & \hbox{for $a\geq 0$}, \\
    \sup_{y\in[a,0]} Z_y, & \hbox{for $a< 0$}, \\
\end{array}%
\right.\vspace{0.2cm}\\
\underline{Z}_a &=& F_a^{(2)}(Z) &=&
\left\{%
\begin{array}{ll}
    \inf_{y\in[0,a]} Z_y, & \hbox{for $a\geq 0$}, \\
    \inf_{y\in[a,0]} Z_y, & \hbox{for $a< 0$}, \\
\end{array}%
\right.\vspace{0.2cm}\\
Z^*_a &=& F_a^{(3)}(Z) &=&
\left\{%
\begin{array}{ll}
    \sup_{y\in[0,a]} |Z_y|, & \hbox{for $a\geq 0$}, \\
    \sup_{y\in[a,0]} |Z_y|, & \hbox{for $a< 0$}, \\
\end{array}%
\right.\vspace{0.2cm}\\
Z^{R}_{a} &=& F_a^{(4)}(Z) &=& Z_a - \underline{Z}_a,\vspace{0.2cm}\\
Z^\#_{a} &=& F_a^{(5)}(Z) &=&
\left\{%
\begin{array}{ll}
    \sup_{0\leq y \leq a}Z^{R}_y, & \hbox{for $a\geq 0$}, \\
    \sup_{a\leq y \leq 0}Z^{R}_y, & \hbox{for $a < 0$}, \\
\end{array}%
\right.\vspace{0.2cm}\\
\sigma_{Z}(a) &=& F_a^{(6)}(Z) &=&
\left\{%
\begin{array}{ll}
\inf\left( x\geq 0 \hbox{ , }Z_{x}\geq a \right), & \hbox{for $a\geq 0$},\\
\inf\left( x\geq 0 \hbox{ , }Z_{x}\leq a \right), & \hbox{for $a< 0$},\\
\end{array}%
\right.\vspace{0.2cm}\\
\widetilde{\sigma}_{Z}(a) &=& F_a^{(7)}(Z) &=& \label{ddd}
\left\{%
\begin{array}{ll}
\inf\left( x\geq 0 \hbox{ , }Z_{-x}\geq a \right), & \hbox{for $a\geq 0$},\\
\inf\left( x\geq 0 \hbox{ , }Z_{-x}\leq a \right), & \hbox{for $a< 0$},\\
\end{array}%
\right.\vspace{0.2cm}\\
U_{Z}(a) &=& F_a^{(8)}(Z) &=& a - \underline{Z}(\sigma_{Z}(a)), \hbox{ for $a\geq 0$},\vspace{0.2cm} \\
\widetilde{U}_{Z}(a) &=& F_a^{(9)}(Z) &=&
a - \underline{Z}\Pare{\widetilde{\sigma}_{Z}(a)} , \hbox{ for $a\geq 0$},\vspace{0.2cm}\\
\widetilde{G}_{Z}(a) &=& F_a^{(10)}(Z) &=&
\widetilde{U}_{Z}(\overline{Z}_{a})\vee Z^\#_a , \hbox{ for $a\geq
0$}.
\end{array}%
$$
Let $\mathcal{D}_i(a)$ for $i\in\{1,\cdots,10\}$ denote the set of
discontinuity points in $\D(\R,\R)$ of $F^{(i)}_a$ and for $v\geq 1$
let $\V^{(v)} = (\V_{vx}/a(v)\hbox{ , }x\in\R)$. From a theorem of
Skorohod \cite{Skoro1}, assumption \ref{hyp1} implies that
$(\V^{(v)}\hbox{ , }v\geq 1)$ converges in law in the Skorohod space
towards $\S$ as $v\to\infty$. It remains to check that the previously
defined functionals have nice continuous properties (with respect
to $\S$) in order to obtain results such as
$F^{(i)}_a(\V^{(v)})\rightarrow F^{(i)}_a(\S)$ in law as
$v\to\infty$.

For $Z\in\D(\R,\R)$ and $a\in\R$, we will say that:
\begin{eqnarray*}
\hbox{$Z$ is oscillating at $a^-$ if }&\forall \varepsilon>0\hbox{  }
&\inf_{(a-\varepsilon,a)}Z < Z_{a-} <\sup_{(a-\varepsilon,a)}Z.\\
\hbox{$Z$ is oscillating at $a^+$ if }&\forall \varepsilon>0\hbox{
} &\inf_{(a,a+\varepsilon)}Z < Z_{a+} <\sup_{(a,a+\varepsilon)}Z.
\end{eqnarray*}
The following lemma collects some easy results about the sample path of $\S$
\begin{lemm} \label{pathS}The following hold:
\begin{enumerate}
\item $\sup_{[0,\infty)}\S = \sup_{(-\infty,0]}\S = \infty$ almost surely.
\item With probability $1$, any path of $\S$ is such that if $\S$ is
discontinuous at a point $x$, then $\S$ is oscillating at $x^-$ and $x^+$.
\item For any fixed $a\in\R$, $\S$ is almost surely continuous at $a$ and
oscillating at $a^-$ and $a^+$.
\end{enumerate}
\end{lemm}
\begin{proof}
(1) and (2) come from Lemma $3.1$ of \cite{Tanak1}, p531 as for
$(3)$, it is well known that $\S$ is almost surely continuous at any
given point and the fact that it is oscillating follows from the
assumption that $|\S|$ is not a subordinator.
\end{proof}
Note that $(2)$ implies that, almost surely, $\S$ is continuous at
all its local extrema. $(2)$ also implies that with
probability $1$, $\S$ attains its bound on any compact interval.
These facts enable us to prove the following:
\begin{prop} \label{fonction1}
For any fixed $a\in\R$ and $i\in\{1,\cdots,10\}$
\begin{equation*}
\Prob{\S \in \mathcal{D}_i(a)} = 0.
\end{equation*}
\end{prop}
\begin{proof}
Let $a$ be fixed. The functionals $F_{i}(a) , i\in\{1,2,3,4,5\}$ are
continuous at all $Z\in\D(\R,\R)$ such that $Z$ is continuous at point $a$
(refer to Proposition 2.11 p305 of \cite{Jacod1} for further
details) and the result follows from (3) of the previous lemma. It is
also easily checked from the definition of the Skorohod topology
that the functionals $F_{i}(a) , i\in\{6,8\}$ are continuous at all
$Z$ which have the following properties:
\begin{enumerate}
\item[(a)] $\sigma_{Z}(a) < \infty$,
\item[(b)] $Z$ is oscillating a $\sigma_Z(a)+$,
\item[(c)] $Z$ attains its bounds on any compact interval.
\end{enumerate}
Using again the previous lemma, we see that (a) and (c) hold for
almost any path of $\S$. Notice that, from the Markov property, part
$(3)$ of the lemma is unchanged when $a$ is replaced by a arbitrary
stopping time hence $(b)$ is also true for almost any path of $\S$.
The proof for $F_{i}(a) , i\in\{7,9\}$ is of course similar.
Finally, the result for $F_{10}(a)$ may easily be deduced from
previous ones using the independence of $(\S_x \hbox{ , }x\geq 0)$
and $(\S_{-x}\hbox{ , }x\geq 0)$.
\end{proof}
We will also use the fact that the random variables $F_{i}(a)$ have
continuous cumulative functions (except for the degenerated cases
$a=0$).
\begin{prop} \label{fonction2}For all $a\neq 0$ and $b\in\R$ and $i\in\{1,\cdots,10\}$:
\begin{equation*}
\Prob{F^{(i)}_a(\S) = b} = 0.
\end{equation*}
\end{prop}
We skip the proof as this may be easily checked from the facts that
$\S$ has a continuous density and the assumption that it is not a
subordinator.

Finally, throughout the rest of this paper, the notation $\Cste{i}$
will always denote a finite strictly positive constant
depending only on our choice of $\P$. In the case of a constant depending on
some other parameters, these will appear in the subscript. We will
also repeatedly use the following lemma easily deduced from the
Uniform Convergence Theorem for regularly varying functions
\cite{Bingh1}, p22 combined  monotonicity property.
\begin{lemm} \label{regvar}
Let $f : [1,\infty)\mapsto \R_{+}$ be a strictly positive non
decreasing function which is regularly varying at infinity with index $\beta \geq
0$. Then, for any $\varepsilon > 0$ there exist
$\Cste{1,\varepsilon, f} , \Cste{2,\varepsilon , f}$ such that for
any $1\leq x \leq y$:
\begin{equation*}
\Cste{1,\varepsilon , f}\Pare{\frac{x}{y}}^{\beta + \varepsilon}
\leq \frac{f(x)}{f(y)}\leq \Cste{2,\varepsilon ,
f}\Pare{\frac{x}{y}}^{\beta - \varepsilon}.
\end{equation*}
\end{lemm}

\subsection{Supremum of the reflected process}
In this subsection, we give some bounds and asymptotics about
$\V^\#$. These estimates which may look quite technical will play a
central role in the proof of Theorem \ref{MainTheo1}. This
subsection is devoted to the proofs of the three following
propositions
\begin{prop} \label{estimonte}We have
\begin{equation*}
\lim_{
\begin{array}{c}
\scriptstyle{x\rightarrow\infty}\\
\scriptstyle{v/a^{-1}(x)\rightarrow\infty}
\end{array}
} \frac{a^{-1}(x)}{v}\log \Prob{\V^{\#}_v \leq x} = -K^\#
\end{equation*}
where $K^\# = -\lim_{v\to\infty } \frac{1}{v}\log \Prob{\S^{\# }_{v}\leq 1}$
is strictly positive and finite.
\end{prop}

\begin{prop} \label{Encadr1}
for all $0< b < 1$, there exists $\Cste{3,b}>0$ such that for all $x$
large enough (depending on $b$) and all $v>0$:
\begin{equation*}
\Cste{3,b}\Prob{\V_v^{\#} \leq x} \leq \Prob{\V_v^{\#} \leq x ,
\overline{\V}_v \leq bx} \leq \Prob{\V_v^{\#} \leq x}.
\end{equation*}
\end{prop}

\begin{prop} \label{Encadr2}
There exists $\Cste{4}>0$ such that for all $x$ large enough and all
$v_1 , v_2 >0$:
\begin{equation*}
\Cste{4}\Prob{\V_{v_1}^{\#} \leq x} \Prob{\V_{v_2}^{\#} \leq x}\leq
\Prob{\V_{v_1 + v_2}^{\#} \leq x}.
\end{equation*}
\end{prop}
Notice that using Proposition \ref{Encadr1} we deduce that
Proposition \ref{estimonte} is unchanged if we replace $\P(\V^{\#}_v
\leq x)$ by $\P(\V^{\#}_v \leq x \hbox{ , }\overline{\V}_v \leq bx)$
for all $b>0$. The proof of the first proposition relies on the
following lemma:
\begin{lemm} \label{lemmlogS} There exists a constant $K^\# \in (0,\infty)$
such that, for any $a,c>0$ and any $b\geq 0$
\begin{equation*}
\lim_{t\to\infty } \frac{a^\alpha}{t}\log \Prob{\S^{\# }_{t}  \leq a
\hbox{ , } \underline{\S}_t \leq -b \hbox{ , } \S_t -
\underline{\S}_t \leq c} = - K^\#.
\end{equation*}
In particular $K^\# =
-\lim_{v\rightarrow\infty}\frac{1}{v}\log(\P(\S^\#_v \leq 1))$.
\end{lemm}
\begin{proof} Using the scaling property, we only need to prove the lemma
in the case  $a=1$. For the sake of clarity, let
\begin{equation*}
\Event{1} = \left\{ \S^{\# }_{t}  \leq 1 \hbox{ , } \underline{\S}_t
\leq -b \hbox{ , } \S_t - \underline{\S}_t \leq c\right\},
\end{equation*}
and let $f(t)= \log \P(\S^{\#}_t \leq 1)$. Using the Markov property
of the stable process $\S$, we deduce that $f(t+s) \leq f(t) + f(s)$
for any $s,t\geq 0$. Since $f$ is subadditive, elementary analysis
shows that the limit $K^\# = -\lim_{t\rightarrow\infty} f(t)/t$
exists and furthermore $K^\# \in (0,\infty]$. In order to prove that
$K^\# <\infty$, note that $\{ \S^{\#}_{t} \leq 1\} \supset
\{\S^*_t \leq 1/2\}$ which implies $f(t)/t \geq
\log \Prob{ \S^*_t\leq 1/2}/t$. Using
Proposition $3$ of \cite{Berto1}, p220, the r.h.s. of this last
inequality converges to some finite constant when $t$ converges to
infinity therefore $K^\#$ must be finite. So we have obtained
\begin{equation*}
\limsup_{t\to\infty } \frac{1}{t}\log \Prob{\Event{1}} \leq
\lim_{t\rightarrow\infty}\frac{1}{t}f\Pare{t}  \leq - K^\#.
\end{equation*}
It remains to prove the lower bound. Let $0 < \varepsilon <\min\Pare{c,1}$ and let $t>1$. Define
\begin{equation*}
\begin{array}{l}
\Event{2} = \left\{ \S_{t-1}^{\#} \leq 1-\varepsilon \right\}, \\
\Event{3} = \left\{ \Pare{\theta_{t-1}\S}^{\#}_{1} \leq \varepsilon
\hbox{ , } \underline{\Pare{\theta_{t-1}\S}}_{1} \leq -b-1 \right\}.
\end{array}
\end{equation*}
We have $\Event{1} \supset \Event{2}\cap \Event{3}$. Since $\S$ has
independent increments, $\Event{2}$ and $\Event{3}$ are independent.
Therefore $\Prob{\Event{1}} \geq \Prob{\Event{2}}\Prob{\Event{3}}$.
Furthermore, using scaling, $\P(\Event{2}) =
f\Pare{(t-1)/(1-\varepsilon)^{\alpha}}$. Hence
\begin{equation} \label{inter6b}
\frac{1}{t}\log \Prob{\Event{1}} \geq \frac{\log\Prob{\Event{3}}}{t} +
\frac{1}{t}f\Pare{\frac{t}{(1-\varepsilon)^{\alpha}}},
\end{equation}
and $\P(\Event{3})= \P(\S^\#_{1} \leq \varepsilon \hbox{ , }
\underline{\S}_{1} \leq -b-1)$ does not depend on $t$ and is not
zero (this is easy to check since $\S$ is not a subordinator).
Taking the limit in (\ref{inter6b}) we conclude that
\begin{equation*}
\liminf_{t\to\infty}\frac{1}{t}\log \Prob{\Event{1}} \geq
\lim_{t\to\infty}\frac{1}{t}f\Pare{\frac{t}{(1-\varepsilon)^{\alpha}}}
= \frac{-K^\infty}{(1-\varepsilon)^\alpha}.
\end{equation*}
\end{proof}

\begin{proof}[Proof of Proposition \ref{estimonte}.]
Let us choose $\varepsilon >0$. The previous lemma combined with
the scaling property of $\S^{\#}$ give
$$K^\# = -\lim_{y\rightarrow\infty}\frac{1}{y^\alpha}\log \Prob{\S^{\#}_1 < \frac{1}{y}}$$
hence we can choose $y_0 > 0$ such that $\log \P(\S^{\#}_1
\leq 1/y_0) \leq - (K^\# - \varepsilon)y_0^{\alpha}$. Combining
results of Proposition \ref{fonction1} and \ref{fonction2} for the
functional $F^{(3)}$ yield:
\begin{equation*}
\lim_{k\rightarrow\infty} \log \Prob{\frac{1}{a(k)}\V^{\#}_k \leq
\frac{1}{y_0}} = \log \Prob{\S^{\#}_1 \leq \frac{1}{y_0}} \leq
-\Pare{K^\# - \varepsilon}y_0^{\alpha}.
\end{equation*}
Therefore, for all $k$ large enough:
\begin{equation} \label{inter7}
\log \Prob{\frac{1}{a(k)}\V^{\#}_k \leq \frac{1}{y_0}} \leq
-\Pare{K^\# - 2\varepsilon}y_0^{\alpha}.
\end{equation}
Let us choose $k = [a^{-1}\Pare{x y_0} ]+1$ thus (\ref{inter7}) holds whenever $x$ is large
enough. Notice the inclusion
\begin{equation*}
\left\{ \V^\#_v \leq x\right\} \subset \bigcap_{n=0}^{[v/k]-1}\left\{(\theta_{nk}\V)^{\#}_{k}\leq x \right\},
\end{equation*}
hence using the independence and stationarity of the increments of
the random walk  at integer times:
\begin{equation} \label{inter8}
\Prob{\V^{\#}_v \leq x} \leq \Pare{\Prob{\V^{\#}_{k} \leq
x}}^{\left[\frac{v}{k}\right]}.
\end{equation}
Since $a(\cdot)$ is nondecreasing, our choice of $k$ implies
$x/a(k) \leq 1/y_0$, therefore:
$$\Prob{\V^{\#}_k \leq x} \leq
\P\Big(\frac{\V^{\#}_k}{a(k)}\leq \frac{1}{y_0}\Big).$$
Combining this inequality with
(\ref{inter7}) and (\ref{inter8}) yields
\begin{equation*}
\log \Prob{\V^{\#}_v \leq x} \leq -\left[\frac{v}{k}\right]
y_0^{\alpha}\Pare{K^\# - 2\varepsilon}.
\end{equation*}
It is easy to check from the regular variation of $a^{-1}(\cdot)$ with
index $\alpha$ that $[v/k]y_0^{\alpha} \sim v/a^{-1}(x)$ when $x$
and $v/a^{-1}(x)$ both go to infinity hence:
\begin{equation*}
\limsup  \frac{a^{-1}(x)}{v}\log \Prob{\V^{\#}_v \leq x} \leq -K^\#.
\end{equation*}
The proof of the lower bound
is quite similar yet slightly more technical. Using Lemma \ref{lemmlogS}
and the scaling property, we can find $y_0> 0$ such that:
\begin{eqnarray}\label{unplk}
\log \Prob{\S^{\#}_{1} \leq \frac{1-\varepsilon}{y_0} \hbox{ , }
\underline{\S}_{1}\leq -\frac{2\varepsilon}{y_0}  \hbox{ , } \S_{1}
- \underline{\S}_{1} \leq \frac{\varepsilon}{y_0}} \geq -\frac{K^\#
y_0^{\alpha}}{(1-2\varepsilon)^{\alpha}}.
\end{eqnarray}
Let us set
\begin{equation*}
\Event{4}\Pare{ k } = \left\{\frac{\V^{\#}_{k}}{a(k)} \leq
\frac{1-\varepsilon}{y_0} \hbox{ , }
\frac{\underline{\V}_{k}}{a(k)}\leq -\frac{2\varepsilon}{y_0} \hbox{
, } \frac{\V_{k} - \underline{\V}_{k}}{a(k)} \leq
\frac{\varepsilon}{y_0}\right\}.
\end{equation*}
Proposition \ref{fonction1} states that the set of continuity points
of the functional:
\begin{eqnarray*}
\D([0,\infty),\R)&\rightarrow&\R^{3}\\
Z&\mapsto&\Pare{Z^\#_1 , \underline{Z}_1 , Z_1 - \underline{Z}_1}
\end{eqnarray*}
has probability $1$  with respect to $\S$. Using Proposition
\ref{fonction2} we deduce
\begin{equation*}
\lim_{k\rightarrow\infty}\Prob{\Event{4}\Pare{ k }} =
\Prob{\S^{\#}_{1} \leq \frac{1-\varepsilon}{y_0} \hbox{ , }
\underline{\S}_{1}\leq -\frac{2\varepsilon}{y_0}  \hbox{ , } \S_{1}
- \underline{\S}_{1} \leq \frac{\varepsilon}{y_0}},
\end{equation*}
hence for all $k$ large enough, it follows from (\ref{unplk}) that
\begin{eqnarray} \label{inter9}
\log \Prob{ \Event{4}\Pare{k}} \geq \frac{-K^\#
y_0^{\alpha}}{(1-3\varepsilon)^{\alpha}}.
\end{eqnarray}
We now choose $k = [a^{-1}\Pare{x y_0}]$. Notice 
that $1/y_0 \leq x/a(k) \leq 2/y_0$ for all $x$ large enough thus
\begin{equation*}
\Event{4}(k) \subset \left\{ \V^{\#}_k \leq (1-\varepsilon)x \hbox{
, }\underline{\V}_{k} \leq -\varepsilon x \hbox{ , } \V_k -
\underline{\V}_k \leq \varepsilon x \right\}.
\end{equation*}
One may check by induction that
\begin{eqnarray*}
\left\{ \V^\#_v \leq x \right\} &\supset&
\bigcap_{n=0}^{[v/k]}\Big\{(\theta_{nk}\V)_k^\# \leq
(1-\varepsilon)x \hbox{ , }\underline{(\theta_{nk}\V)}_k \leq -
\varepsilon x , \\
& &\hspace{1cm} (\theta_{nk}\V)_k -
\underline{(\theta_{nk}\V)}_k\leq \varepsilon x\Big\},
\end{eqnarray*}
hence using independence and stationarity of the increments of $\V$
a integer times:
\begin{eqnarray*}
\Prob{\V^{\#}_v \leq x} &\geq& \P\Big( \V^{\#}_k \leq
(1-\varepsilon)x \hbox{ , } \underline{\V}_{k} \leq -\varepsilon x
\hbox{ , } \V_k -
\underline{\V}_k \leq \varepsilon x \Big)^{\left[\frac{v}{k}\right]+1} \\
& \geq & \Prob{\Event{4}(k)}^{\left[\frac{v}{k}\right]+1}.
\end{eqnarray*}
Combining this inequality with (\ref{inter9}) we get for any $x$ large
enough:
\begin{equation*}
\log \Prob{\V^{\#}_v \leq x } \geq
\frac{-K^\#}{(1-3\varepsilon)^{\alpha}}\left(\left[\frac{v}{k}\right]+1\right)y_0^\alpha.
\end{equation*}
Notice that $([v/k]+1)y_0^\alpha \sim v/a^{-1}(x)$ as $x$ and
$v/a^{-1}\Pare{x}$ go to infinity simultaneously which completes the
proof.
\end{proof}

\begin{proof}[Proof of Proposition \ref{Encadr1}.]
The upper bound is trivial. Let $0<b<1$, define $v_1 = [a^{-1}(x)]$ and set $c = (b-1)x$:
\begin{eqnarray*}
\left\{ \V^\#_v \leq x , \overline{\V}_v \leq bx \right\} &\supset&
\left\{ \V^\#_v \leq x , \overline{\V}_v \leq bx , \sigma_{\V}(c)
\leq v_1 \right\} \\
&\supset& \left\{ \V^\#_{\sigma_{\V}(c)} \leq bx , \sigma_{\V}(c) \leq v_1
\right\} \cap \left\{ \Pare{\theta_{\sigma_{\V}(c)}\V}^\#_{v}
\leq x\right\},
\end{eqnarray*}
thus
\begin{eqnarray*}
\Prob{\V^\#_v \leq x , \overline{\V}_v \leq bx} &\geq&
\Prob{\V^\#_{\sigma_{\V}(c)} \leq bx , \sigma_{\V}(c) \leq v_1}
\Prob{\V^\#_v \leq x}\\
& \geq & \Prob{\V^\#_{v_1} \leq bx ,\underline{\V}_{v_1} \leq c} \Prob{\V^\#_v \leq x}.
\end{eqnarray*}
Just like for the previous proof, we see that $\P(\V^\#_{v_1} \leq bx
, \underline{\V}_{v_1} \leq c)$
converges when $x$ goes to infinity toward $\P(\S^\#_{1} \leq b ,
\underline{\S}_{1} \leq b-1)$ and this quantity is
strictly positive number because $|\S|$ is not a subordinator.
\end{proof}
\begin{proof}[Proof of Proposition \ref{Encadr2}.] Notice that
\begin{eqnarray*}
\left\{\V_{v_1 + v_2}^{\#} \leq x\right\} &\supset& \left\{
\V^\#_{[v_1]+[v_2]+2}\leq x\right\}\\
&\supset& \Big\{\V_1 \leq 0 \hbox{ , } \V_2 - \V_1 \leq
0\Big\}\\
&&\cap\left\{(\theta_{2}\V)^\#_{[v_1]}\leq x \hbox{ , }
(\theta_{2}\V)_{[v_1]}-\underline{(\theta_{2}\V)}_{[v_1]} \leq \frac{x}{2}\right\}\\
&&\cap\left\{ (\theta_{2+[v_1]}\V)^\#_{{[v_2]}}\leq x\hbox{ ,
}\overline{(\theta_{2+[v_1]}\V)}_{[v_2]}\leq\frac{x}{2} \right\},
\end{eqnarray*}
Using the independence and stationarity of the increments of $\V$ at integer
time and setting $\Cste{5} = \P(\V_1 \leq 0)>0$ we see
that $\P(\V_{v_1 + v_2}^{\#} \leq x)$ is greater than
\begin{equation*}
\Cste{5}^2\Prob{
\V^\#_{[v_1]}\leq x\hbox{ , }\V_{[v_1]}-\underline{\V}_{[v_1]}\leq
\frac{x}{2}}\Prob{ \V^\#_{[v_2]}\leq x\hbox{ ,
}\overline{\V}_{[v_2]}\leq \frac{x}{2}},
\end{equation*}
but time reversal of the random walk $\V$ shows that:
$$\Prob{\V^\#_{[v_1]}\leq x , \V_{[v_1]} - \underline{\V}_{[v_1]} \leq
x/2} = \Prob{\V^\#_{[v_1]}\leq x , \overline{\V}_{[v_1]}\leq x/2},$$
hence using Proposition \ref{Encadr1}:
\begin{eqnarray*}
\P(\V_{v_1 + v_2}^{\#} \leq x) &\geq& (\Cste{3,\frac{1}{2}}
\Cste{5})^2\Prob{ \V^\#_{[v_1]}\leq x}\Prob{ \V^\#_{[v_2]}\leq x} \\
&\geq& (\Cste{3,\frac{1}{2}}
\Cste{5})^2\Prob{ \V^\#_{v_1}\leq x}\Prob{ \V^\#_{v_2}\leq x}.
\end{eqnarray*}
\end{proof}

\subsection{The case where $\S$ is a completely asymmetric stable process.}\label{subsectionA}
One may wish to calculate the value of the constant $K^\#$ that
appears in the last section. Unfortunately, we do not know its value in
general. However, the completely asymmetric case is a particularly
nice setting where calculations may be carried to their full extend. We
now assume throughout this section that the stable process
$(\S_x\hbox{ , }x\geq 0)$ either has no positive jumps
hence the exponential moments of $\S$ are finite and (\ref{ee1})
hold (recall that we assume $\gamma'=1$) or $\S$ has no negative
jumps thus $\Esp{\exp(-\lambda\S_t)} = \exp(t\lambda^\alpha)$ for all
$t,\lambda\geq 0$. For $a,b>0$, define the stopping times:
\begin{eqnarray*}
\tau_b &=& \inf(t\geq 0\hbox{ , }\S_t \geq b) = \sigma_{\S}(b),\\
\tau^\#_b &=& \inf(t\geq 0\hbox{ , } \S^\#_t \geq b) =
\sigma_{\S^\#}(b),\\
\tau^*_{a,b} &=& \inf (t\geq 0\hbox{ , } \S_t \hbox{ not in }
(-a,b)).
\end{eqnarray*}
Recall that $\Mittag_\alpha$ stand for the Mittag Leffler function with
parameter $\alpha$.
\begin{prop} \label{laplaceasy}
When $\S$ has no positive jumps:
\begin{equation*}
\Esp{e^{-q\tau^\#_1}} = \frac{1}{\Mittag_\alpha(q)},
\end{equation*}
and when $\S$ has no negative jumps:
\begin{equation*}
\Esp{e^{-q\tau^{\#}_1}} = \Mittag_\alpha(q) - \frac{\alpha q
(\Mittag_\alpha '(q))^2}{\alpha q \Mittag_\alpha''(q) + (\alpha
-1)\Mittag_\alpha'(q)}.
\end{equation*}
\end{prop}
This proposition is a particular case from Proposition $2$ of
\cite{Pisto1} p191. Still, we give here a simpler proof when $\S$ is
stable using the solution of the two sided exit problem given by
Bertoin in \cite{Berto2}.
\begin{proof} We suppose that $\S$ has no negative jumps.
Let $\eta(q)$ be an exponential random time of parameter $q$ independent of
$\S$. Let also $a,b$ be strictly positive real numbers such that
$a+b=1$. We may without loss of generality assume any path of $\S$
attains its bounds on any compact interval and is continuous at all
local extrema (because this happens with probability $1$ according
to Lemma \ref{pathS}) thus on the one hand, the event $\{\tau^{\#}_1 > \eta(q)\}$ contains
\begin{equation*}
\Big\{\tau^*_{a,b} >
\eta(q)\Big\}\cup\Pare{\left\{\tau^*_{a,b}\leq \eta(q)\hbox{ , }
\S_{\tau^*_{a,b}} \leq
-a\right\}\cap\left\{(\theta_{\tau^*_{a,b}}\S)^\#_{\eta(q)-\tau^*_{a,b}}
<1 \right\}}.
\end{equation*}
Using the strong Markov property of $\S$, the lack of memory and
the independence of the exponential time, it follows that $\P(\tau^\#_1 > \eta(q))$ is greater than
\begin{equation*}
\Prob{\tau^*_{a,b} >
\eta(q)}+\Prob{\tau^*_{a,b}\leq \eta(q)\hbox{ ,
}\S_{\tau^*_{a,b}}\leq -a}\Prob{\tau^\#_1 > \eta(q)},
\end{equation*}
therefore
\begin{equation} \label{mont1}
\Prob{\tau^{\#}_1 > \eta(q)} \geq \frac{\Prob{\tau^*_{a,b} >
\eta(q)}}{1-\Prob{\tau^*_{a,b}\leq \eta(q)\hbox{ ,
}\S_{\tau^*_{a,b}}\leq -a}}.
\end{equation}
On the other hand, one may check that the event $\{\tau^{\#}_1 > \eta(q)\}$ is a subset of
\begin{equation*}
\Big\{\tau^*_{a,b} >
\eta(q)\Big\}\cup\Pare{\left\{\tau^*_{a,b}\leq \eta(q)\hbox{ ,
}\S_{\tau^*_{a,b}} \leq
-a\right\}\cap\left\{(\theta_{\tau^*_{a,b}}\S)^\#_{\eta(q)-\tau^*_{a,b}}
<b \right\}},
\end{equation*}
and similarly we deduce
\begin{equation}\label{mont2}
\Prob{\tau^\#_b > \eta(q)} \leq \frac{\Prob{\tau^*_{a,b} >
\eta(q)}}{1-\Prob{\tau^*_{a,b}\leq \eta(q)\hbox{ ,
}\S_{\tau^*_{a,b}}\leq -a}}.
\end{equation}
Obviously $\tau^\#_b$ converges to
$\tau^\#_1$ almost surely as $b$ converges to $1$. Combining
this observation with (\ref{mont1}) and (\ref{mont2}), we find:
\begin{equation} \label{mont3}
\Prob{\tau^\#_1 > \eta(q)} = \lim_{b\nearrow 1}
\frac{\Prob{\tau^*_{1-b,b}
> \eta(q)}}{1-\Prob{\tau^*_{1-b,b}\leq \eta(q)\hbox{ , }\S_{\tau^*_{1-b,b}}\leq b-1}}.
\end{equation}
The value of the probabilities of the r.h.s. of this equation have been
calculated by Bertoin in \cite{Berto2}:
\begin{eqnarray}
\label{Bertcal1}&\Prob{\tau^*_{1-b,b} > \eta(q)} = 1 -
\Mittag_\alpha(b^{\alpha}) +
\frac{b^{\alpha-1} \Mittag_\alpha'(qb^\alpha)}{\Mittag_\alpha'(q)}\Pare{\Mittag_\alpha(q)-1},& \\
\label{Bertcal2}&\Prob{\tau^*_{1-b,b}\leq \eta(q)\hbox{ ,
}\S_{\tau^*_{1-b,b}}\leq b-1} =
\frac{b^{\alpha-1}\Mittag'_\alpha(qb^\alpha)}{\Mittag'_\alpha(q)}.&
\end{eqnarray}
A Taylor expansion of $\Mittag_\alpha$ and $\Mittag_\alpha'$ near
point $q$ enables us to calculate the limit in (\ref{mont3}) in term
of $\Mittag_\alpha$ and its first and second derivatives. After a few
lines of elementary calculus:
\begin{equation*}
\Prob{\tau^\#_1 > \eta(q)} = 1 -\Mittag_\alpha (q) + \frac{\alpha q
(\Mittag_\alpha ' (q))^2}{\alpha q \Mittag_\alpha'' (q) + (\alpha
-1) \Mittag_\alpha' (q)}.
\end{equation*}
We complete the proof using the well known relation
$\EE(\exp(-q\tau^\#_1))= 1 -\P(\tau^\#_1 > \eta(q))$. The proof in
the case where $\S$ has no positive jumps is similar (and the
calculation of the limit is even easier). We omit it.
\end{proof}
\begin{corr} \label{valeurKdiese}
Recall that $-\rho_1(\alpha)$ is the first negative root of
$\Mittag_\alpha$ and $-\rho_2(\alpha)$ is the first negative root of
$\alpha x \Mittag_\alpha''(x) + (\alpha-1)\Mittag_\alpha'(x)$. The
constant of Proposition \ref{estimonte} is given by:
\begin{equation*}
K^\# =
\left\{%
\begin{array}{l}
    \rho_1 (\alpha) \hbox{ when $\S$ has no positive jumps}, \\
    \rho_2 (\alpha) \hbox{ when $\S$ has no negative jumps}. \\
\end{array}%
\right.
\end{equation*}
\end{corr}
\begin{proof} Recall that $K^\# = -\lim_{t\rightarrow\infty}\P(\S^\#_t \leq 1)/t$.
Using the same argument as in Corollary $1$ of \cite{Berto2}, we see
that $K^\# = \rho_1(\alpha)$ when $\S$ has no positive jumps.
Similarly, when $\S$ has no negative jumps $-K^\#$ is equal to the
first negative pole of
\begin{equation*}
g(x) = \frac{\alpha x (\Mittag_\alpha'(x))^2}{\alpha x
\Mittag_\alpha ''(x) + (\alpha -1)\Mittag_\alpha '(x)} =
\Mittag_\alpha (x) - \Esp{e^{-x \tau^\#_1}}.
\end{equation*}
Let $-x_0$ be the first negative root of $\Mittag_\alpha '$. Since
$\Mittag_\alpha '(0)>0$, this implies that $\Mittag_\alpha$ is
strictly increasing on $[-x_0 , 0]$. Note also that $x\mapsto
-\EE(\exp(-x \tau^\#_1 ))$ is increasing on $(-K^\#,0]$ thus $g(x)$
is strictly increasing on $(-(K^\# \wedge x_0) , 0]$. Since
$g(-x_0)=g(0)=0$ (this holds even when $-x_0$ is a zero of multiple
order) we deduce from the monotonicity of $g$ that $K^\# < x_0$ and
this shows that the first negative pole of $g$ is indeed $-\rho_2
(\alpha)$.
\end{proof}
We conclude this subsection by calculating the Laplace transform of
$\tau^\#_{1}\wedge\tau_b$. This will be useful for the determination
of the limiting law in the proof of Theorem \ref{MainTheo3}.
\begin{corr}\label{pourloi}
for $0 < b \leq 1$, when $\S$ has no positive jumps
\begin{equation*}
\Esp{e^{-q \tau^\#_{1}\wedge\tau_b}}
=\frac{\Mittag_\alpha(q(1-b)^{\alpha})}{\Mittag_\alpha(q)},
\end{equation*}
and when $\S$ has no negative jumps
\begin{equation*}
\Esp{e^{-q \tau^\#_{1}\wedge\tau_b}} =\Mittag_\alpha(q b^\alpha) -
b^{\alpha-1}\frac{\alpha q
\Mittag_\alpha'(qb^\alpha)\Mittag_\alpha'(q)}{\alpha q
\Mittag_\alpha''(q) + (\alpha-1)\Mittag_\alpha'(q)}.
\end{equation*}
\end{corr}
\begin{proof} Let $\eta(q)$ still denote an exponential time with parameter $q$
independent of $\S$. Suppose that $\S$ has no negative jumps, using
the Markov property and the lack of memory of the exponential law we
get
\begin{eqnarray*}
\Prob{ \tau^\#_{1}\wedge\tau_b> \eta(q)} &=& \Prob{\tau^{*}_{1-b,b} \leq \eta(q) , \S_{\tau^{*}_{1-b,b}}
\leq b-1}\Prob{\tau^{\#}_1 > \eta(q)}\\
&&+\Prob{\tau^{*}_{1-b,b} > \eta(q)}.
\end{eqnarray*}
The r.h.s. of the last equality may be calculated explicitly using
again (\ref{Bertcal1}), (\ref{Bertcal2}) and Proposition
\ref{laplaceasy} hence, after simplication:
\begin{equation*}
\Prob{\tau_1^\#\wedge\tau_b > \eta(q)} = 1 - \Mittag_\alpha(q
b^\alpha) + b^{\alpha-1}\frac{\alpha q
\Mittag_\alpha'(qb^\alpha)\Mittag_\alpha'(q)}{\alpha q
\Mittag_\alpha''(q) + (\alpha-1)\Mittag_\alpha'(q)}.
\end{equation*}
The no positive jumps case may be treated the same way.
\end{proof}

\subsection{The exit problem for the random walk $\V$}
Let us define for $x,y>0$ the following events:
\begin{eqnarray}
\nonumber\Lambda\Pare{x,y} &= &\left\{ \Pare{\V_s}_{s\geq 0} \hbox{
hits }
(y,\infty) \hbox{ before it hits }(-\infty,-x)\right\}, \\
\nonumber\Lambda'\Pare{x,y} &= &\left\{ \Pare{\V_s}_{s\geq 0} \hbox{
hits }
[y,\infty) \hbox{ before it hits }(-\infty,-x]\right\}, \\
\nonumber\widetilde{\Lambda}'\Pare{x,y} &= &\left\{
\Pare{\V_{-s}}_{s\geq 0} \hbox{ hits } (-\infty,-y] \hbox{ before it
hits }[x,\infty)\right\}.
\end{eqnarray}
We are interested in the behavior of the probabilities of these
events for large $x,y$. In the case of a fixed $x$, when $y$ goes to
infinity, this study was done by Bertoin and Doney in \cite{Doney2}.
Here, we need to study this quantities when both $x$ and $y$ go to
infinity with the ratio $y/x$ also going to infinity. We already
defined $(\T_n)_{n\geq 0}$ to be the sequence of strict descending
ladder times, we now consider the associated ladder heights
$(\H_n)_{n\geq 0}$:
\begin{equation*}
\H_{n} = -\V_{\T_n}.
\end{equation*}
We will also need the sequence $(\M_{n})_{n\geq 1}$:
\begin{equation*}
\M_n = \max\Pare{\V_{k} + \H_{n-1} \hbox{ , } \T_{n-1}\leq k <
\T_{n}}.
\end{equation*}
Note that the sequence $\Pare{\T_{n+1}-\T_n, \H_{n+1} - \H_n ,
\M_n}_{n\geq 1}$ is independent, identically distributed. We know
that $\T_1$ is in the domain of attraction of a positive stable law
of index $q$ with norming constants $(b_n)$. Now Corollary $3$ of
\cite{Doney1} gives $\Prob{\M_1 > x}$ regularly varying with index
$-\alpha q$. More precisely, it gives:
\begin{equation} \label{tail_M1}
\Prob{\M_1 > x} \underset{x\to\infty}{\sim} \frac{\Cste{6}}{
b^{-1}\Pare{a^{-1}\Pare{x}}}.
\end{equation}
In particular, this shows that $\M_1$ is in the domain of attraction
of a positive stable law when $\alpha q <1$ and that $\M_1$ is
relatively stable when $\alpha q =1$ (relatively stable means that
$\frac{1}{a(b(n))}\sum_{k\leq n} \M_k$ converges in probability to
some strictly positive constant).

For $\H_1$, using Theorem $9$ of \cite{Rogoz1}, we see that $\H_1$
is in the domain of attraction of a positive stable law with index
$\alpha q$ when $\alpha q <1$ and that $\H_1$ is relatively stable
in the case $\alpha q =1$. Furthermore, the lemma of \cite{Doney1},
p358 shows that we can choose $a\Pare{b\Pare{n}}$ as norming
constant for $\H_1$ in any of those two cases. That is:
\begin{equation*}
\frac{\H_n}{a\Pare{b\Pare{n}}} \hbox{ converges to }
\left\{%
\begin{array}{l}
\hbox{some constant $\Cste{7}$, in probability when $\alpha q=1$,} \\
\hbox{a positive stable law of index $\alpha q$ otherwise.}
\end{array}%
\right.
\end{equation*}
When $\alpha q <1$, this shows that (\ref{tail_M1}) holds with
$\H_1$ in place of $\M_1$ (for a different value of $\Cste{6}$).
Unfortunately, in the case $\alpha q = 1$, the relative stability of
$\H_1$ does not imply the regular variation of $\Prob{\H_1 >x}$
(look at the counter example in \cite{Rogoz1}, p 576). However, we
can still prove a smooth behavior for the associated renewal
function:
\begin{equation*}
\RenH(x) = \sum_{n=0}^{\infty} \Prob{\H_n \leq x}.
\end{equation*}
\begin{lemm} \label{tail_U}
there exists a constant $\Cste{8} >0$ such that
\begin{eqnarray*}
\RenH\Pare{x}\underset{x\to\infty}{\sim} \Cste{8}
b^{-1}\Pare{a^{-1}\Pare{x}}.
\end{eqnarray*}
\end{lemm}
\begin{proof} When $\alpha q <1$ we mentioned that $\Prob{\H_1 > x} \sim
\Cste{9}/b^{-1}\Pare{a^{-1}\Pare{x}}$ where $\Cste{9}$ is some
strictly positive constant. In this case, the asymptotic behavior of
$\RenH$ follows from the Tauberian Theorem as in Lemma p446 of
\cite{Felle1}. We now consider the case $\alpha q = 1$. Let
$L(\lambda) = \Esp{e^{-\lambda \H_1}}$ stand for the Laplace
transform of $\H_1$. We know that
\begin{equation*}
\frac{\H_n}{a\Pare{b\Pare{n}}}
\underset{n\to\infty}{\overset{\hbox{Prob.}}{\longrightarrow}} \Cste{7}
\end{equation*}
therefore, for any $\lambda \geq 0$ and when $n$ ranges trough the
set of integers:
\begin{equation} \label{inter1}
\left(L\Pare{\frac{\lambda}{a\Pare{b\Pare{n}}}}\right)^{n}
\underset{n\to\infty}{\longrightarrow} e^{-\Cste{7}\lambda}.
\end{equation}
Since $L$ is continuous at $0$ with $L(0)=1$, setting $\lambda =
1$ and taking the logarithm in (\ref{inter1}) give
\begin{equation} \label{inter2}
n\Pare{1 - L\Pare{\frac{1}{a\Pare{b\Pare{n}}}}}
\underset{n\to\infty}{\longrightarrow} \Cste{7}.
\end{equation}
Using the monotonicity of $L$ and $a\Pare{b\Pare{\cdot}}$, it is
easy to check that (\ref{inter2}) still holds when $n$ now ranges
trough the set of real numbers, thus:
\begin{equation} \label{inter3}
1 - L\Pare{\frac{1}{x}}\underset{x\to\infty}{\sim}\frac{\Cste{7}}{b^{-1}\Pare{a^{-1}\Pare{x}}}.
\end{equation}
Let us now define $\widehat{\RenH}(y) = \int_{0}^{\infty} e^{-y
x}\RenH\Pare{dx}$. The well-known relation $\widehat{\RenH}\Pare{y}
= 1/\Pare{1 - L(y)}$ combined with (\ref{inter3}) shows that
$\widehat{\RenH}$ is regularly varying near $0$ hence we can use
Karamata's Tauberian/Abelian Theorem to conclude the proof.
\end{proof}

\begin{prop} \label{tempssortie} There exists $\Cste{10}$ such that
when $x\to\infty$ and $\frac{y}{x}\to\infty$,
\begin{equation*}
\Prob{\Lambda\Pare{x,y}} \sim \Cste{10}
\frac{b^{-1}\Pare{a^{-1}\Pare{x}}}{b^{-1}\Pare{a^{-1}\Pare{x+y}}}.
\end{equation*}
This result also hold for  $\P(\Lambda'\Pare{x,y})$ and
$\P(\widetilde{\Lambda}'\Pare{x,y})$.
\end{prop}

\begin{proof} The two processes $(\V_s)_{s\geq 0}$ and $(-\V_{-s})_{s\geq 0}$ have
the same law hence $\P(\Lambda'\Pare{x,y})= \P(\widetilde{\Lambda}'\Pare{x,y})$. We
also have the trivial inclusion $\Lambda\Pare{x-1 , y} \subset
\Lambda'\Pare{x,y} \subset \Lambda\Pare{x,y-1}$, so we only need to
prove the proposition for $\Lambda\Pare{x,y}$. The first part of the
proof is borrowed from  Bertoin and Doney \cite{Doney2}, p2157. The probability 
$\P(\Lambda(x,y))$ is equal to
\begin{eqnarray}\label{unrefdeplus}
\nonumber\Prob{\M_1 > y} &+& \sum_{k=1}^{\infty}
\P\Big(\M_1 \leq y+ \H_0 , \cdots , \M_k \leq y + \H_{k-1},\\
&&\hspace{1.2cm} \H_k \leq x , \M_{k+1} > y + \H_k\Big),
\end{eqnarray}
thus
\begin{eqnarray*}
\Prob{\Lambda\Pare{x,y}}  & \leq  & \Prob{\M_1 > y} + \sum_{k=1}^{\infty}\Prob{\H_k \leq x , \M_{k+1} > y + \H_k} \\
& \leq & \Prob{\M_1 > y} + \sum_{k=1}^{\infty}\Prob{\H_k \leq x , \M_{k+1} > y} \\
& \leq & \Prob{\M_1 > y}\RenH(x).
\end{eqnarray*}
Using (\ref{tail_M1}), Lemma \ref{tail_U} and the equivalence
$\P(\M_1 > y)\sim\P(\M_1> x+y)$ when $x$ and $y/x$ go to infinity,
we obtain the upper bound with $\Cste{10} = \Cste{6} \Cste{8}$.
We now prove the result pertaining to the lower bound. Let $k_0\in \N^{*}$.
From (\ref{unrefdeplus}), we see that $\P(\Lambda(x,y))$ is bigger than
\begin{eqnarray*}
&&\Prob{\M_1 > y} + \sum_{k=1}^{\infty}\Prob{\M_1 \leq y ,\cdots , \M_k \leq y , \H_k \leq x , \M_{k+1} > x+y} \\
&&\qquad\geq \Prob{\M_1 > x+ y}\Big(1 + \sum_{k=1}^{k_0} \Prob{\M_1 \leq y , \cdots , \M_k \leq y , \H_k \leq x}\Big),
\end{eqnarray*}
hence
\begin{equation}\label{inter4}
\Prob{\Lambda\Pare{x,y}} \geq \Prob{\M_1 > x+y}\Big(\RenH\Pare{x} - \RenH_{k_0}\Pare{x} -
\WartW_{k_0}\Pare{y}\Big),
\end{equation}
with
\begin{eqnarray*}
\RenH_{k_0}\Pare{x}  &=& \sum_{k= k_0 +1}^{\infty}\Prob{\H_k \leq x}, \\
\WartW_{k_0}\Pare{y}  &=& \sum_{k=1}^{k_0}\Prob{\M_1 > y \hbox { or
} \cdots \hbox{ or } \M_{k} > y }.
\end{eqnarray*}
On the one hand, using (\ref{tail_M1}) and Lemma \ref{tail_U}, for $y$
large enough:
\begin{equation*}
\WartW_{k_0}\Pare{y}  \leq  \sum_{k=1}^{K}k_0 \Prob{\M_1 > y} \leq
k_0^2 \Prob{\M_1 > y} \leq \frac{\Cste{11} k_0^2}{\RenH\Pare{y}}.
\end{equation*}
On the other hand:
\begin{eqnarray*}
\RenH_{k_0}\Pare{x} & = & \sum_{k=0}^{\infty}\Prob{\H_{k_0+1} + \Pare{\H_{k+k_0+1} - \H_{k_0+1} \leq x}} \\
& \leq & \sum_{k=0}^{\infty}\Prob{\H_{k+k_0+1} - \H_{k_0+1} \leq x}\Prob{\H_{k_0+1} \leq x}\\
& \leq & \RenH\Pare{x}\Prob{\H_{k_0} \leq x}.
\end{eqnarray*}
Combining these two bounds with (\ref{inter4}) yields, for all $x,y$
large enough:
\begin{equation*}
\Prob{\Lambda\Pare{x,y}}  \geq   \Prob{\M_1 >
x+y}\RenH\Pare{x}\Pare{1  - \Prob{\H_{k_0} \leq x} -\frac{\Cste{11}
k_0^2}{\Pare{\RenH\Pare{y}}^2} }.
\end{equation*}
It only remains to show that for a good choice of $k_0=k_0 (x,y)$, we have
\begin{equation*}
\Prob{\H_{k_0} \leq x} + \frac{\Cste{11}
k_0^2}{\Pare{\RenH\Pare{y}}^2}
\underset{x,\frac{y}{x}\to\infty}{\longrightarrow} 0.
\end{equation*}
Let  $k_0 = \left[b^{-1}\Pare{a^{-1}\Pare{x
\log\Pare{y/x}}}\right]$. Note that $k_0$ is such that
$k_0\to\infty$, when $x$ and
$y/x$ go to infinity simultaneously, and  we know that
\begin{equation*}
\frac{\H_{k_0}}{a\Pare{b\Pare{k_0}}}\underset{k_0\to\infty}{\overset{\hbox{law}}{\longrightarrow}}
J_\infty
\end{equation*}
where $J_\infty$ is either a positive stable law ($\alpha q <1$) or
a strictly positive constant ($\alpha q= 1$). In either cases
$\Prob{J_\infty =0} = 0$. Since $x/a(b(k_0))\to 0$ when $x$ and
$y/x$ go to infinity simultaneously we deduce:
\begin{equation*}
\Prob{\H_{k_0} \leq x} = \Prob{\frac{\H_{k_0}}{a\Pare{b\Pare{k_0}}}
\leq
\frac{x}{a\Pare{b\Pare{k_0}}}}\underset{x,\frac{y}{x}\to\infty}{\longrightarrow}
0.
\end{equation*}
Finally, using Lemmas \ref{regvar} and \ref{tail_U} we conclude that
\begin{equation*}
\frac{\Cste{11} k_0^2}{\Pare{\RenH\Pare{y}}^2}
\underset{x,\frac{y}{x}\to\infty}{\sim}
\frac{\Cste{11}}{\Cste{8}^2}\Pare{\frac{\RenH\Pare{x\log\frac{y}{x}}}{\RenH\Pare{y}}}^{2}
\underset{x,\frac{y}{x}\to\infty}{\longrightarrow} 0.
\end{equation*}
\end{proof}

\subsection{Other estimates}
We conclude the section about the fluctuations of $\V$ by collecting
several results on the functional
$\overline{\V}$ and $\underline{\V}$. We start with a reflection
principle for $\V$:
\begin{lemm} \label{lemmreflex} There exists $\Cste{12}$ such that for all $v,x >0$:
$$\Prob{\overline{\V}_{v} \geq x}\leq \Cste{12} \Prob{\V_{v} \geq x},$$
similarly
$$\Prob{\underline{\V}_{v} \leq -x}\leq \Cste{12} \Prob{\V_{v} \leq -x}.$$
\end{lemm}

\begin{proof}
We only need to prove the first inequality as the second can be
obtained in the same way (with a possibly extended value for
$\Cste{12}$).
\begin{eqnarray*}
\Prob{\overline{\V}_v \geq x} & = & \Prob{\sigma_{\V}(x) \leq [v]} \\
&\leq& \Prob{\sigma_{\V}(x) \leq [v] \hbox{ , }\V_{[v]} < x}+ \Prob{\V_v \geq x} \\
&\leq& \sum_{k=1}^{[v]}\Prob{\sigma_{\V}(x) = k \hbox{ , }\V_{[v]} <
x} + \Prob{\V_v \geq x}.
\end{eqnarray*}
From the Markov property, we check that $\P(\sigma_{\V}(x)=k ,
\V_{[v]}<x)$ is equal to
\begin{eqnarray*}
&&\Prob{\sigma_{\V}(x)=k}\int_{\petit{y\geq x}}
\Prob{\V_{[v]-k}<x-y}\Prob{\V_{\sigma_{\V}(x)}=dy |
\sigma_{\V}(x)=k}\\
&&\qquad\leq\Prob{\sigma_{\V}(x) =
k}\Prob{\V_{[v]-k} <0}.
\end{eqnarray*}
Our assumption on $\V$ implies that
$\lim_{n\rightarrow\infty}\Prob{\V_n < 0} = \Prob{\S <0} = q < 1$
thus, there exists $\Cste{13} >0$ such that $\sup_{n} \Prob{\V_n < 0}
= \Cste{13} < 1$. Therefore
\begin{eqnarray*}
\Prob{\overline{\V}_v \geq x} &\leq& \Cste{13} \sum_{k=1}^{[v]}\Prob{\sigma_{\V}(x) = k } + \Prob{\V_v \geq x} \\
&\leq& \Cste{13}\Prob{\sigma_{\V}(x) \leq v} + \Prob{\V_v \geq x} \\
&\leq& \frac{1}{1-\Cste{13}}\Prob{\V_v \geq x}.
\end{eqnarray*}
\end{proof}
We now estimate the large deviations of $\Prob{\V_v > x}$. Using the
characterization of the domains of attraction to a stable law (see
chapter IX, section 8 of \cite{Felle1}), assumption \ref{hyp1}
implies:
\begin{equation} \label{gd1}
a^{-1}(x)\Prob{\V_1 > x} \underset{x\to\infty}{\longrightarrow}
\left\{%
\begin{array}{ll}
    \Cste{14}>0 & \hbox{ if $\S$ has positive jumps,} \\
    0 & \hbox{ otherwise.}
\end{array}%
\right.
\end{equation}
Similarly:
\begin{equation} \label{gd2}
a^{-1}(x)\Prob{\V_1 < - x} \underset{x\to\infty}{\longrightarrow}
\left\{%
\begin{array}{ll}
    \Cste{15}>0 & \hbox{ if $\S$ has negative jumps,} \\
    0 & \hbox{ otherwise.}
\end{array}%
\right.
\end{equation}
\begin{prop}\label{precislemmefluc}
there exists $\Cste{16}>0$ such that for all $v\geq 1$ and all $x\geq
1$:
\begin{equation}\label{regd1}
\Prob{\V_v > x} \leq \Cste{16}\frac{v}{a^{-1}(x)}.
\end{equation}
Moreover, if $\S$ has positive jumps:
\begin{equation}\label{regd2}
\Prob{\V_v > x} \underset{ \petit{
\begin{array}{c}
v\to\infty \\
\frac{a^{-1}(x)}{v}\to\infty
\end{array}%
}}{\sim} v \Prob{\V_1 > x} \underset{ \petit{
\begin{array}{c}
v\to\infty \\
\frac{a^{-1}(x)}{v}\to\infty
\end{array}%
}}{\sim} \Cste{14}\frac{v}{a^{-1}(x)}.
\end{equation}
There is of course a similar result for $\Prob{\V_v < -x}$.
\end{prop}
\begin{proof}
Result (\ref{regd2}) is already known and is stated in \cite{Borov1}
yet we could not find a proof of this result in English. A weaker
result is proved by Heyde \cite{Heyde1} but a slight modification of his argument
will enable us to prove the proposition. Let us choose $1/2 < \delta < 1$ and set  $z = (x/a(v))^\delta
a(v)$. Define for $ k\geq 1$:
\begin{equation*}
\zeta_{k,z} =
\left\{%
\begin{array}{ll}
    \V_{k}-\V_{k-1} & \hbox{ if $|\V_{k}-\V_{k-1}| \leq z$,} \\
    0 & \hbox{ otherwise.}
\end{array}%
\right.
\end{equation*}
Let $\varepsilon >0$ and set:
\begin{eqnarray*}
\Event{5} &=& \Big\{ \V_{k}-\V_{k-1} > (1-\varepsilon)x\hbox { for at least one $k$ in $\{1,\ldots,[v]\}$}\Big\}, \\
\Event{6} &=& \Big\{ \V_{k}-\V_{k-1} > z\hbox { for at least two $k$'s in $\{1,\ldots,[v]\}$}\Big\}, \\
\Event{7} &=& \Big\{ \zeta_{1,z} + \ldots + \zeta_{[v],z} >
\varepsilon x \Big\}.
\end{eqnarray*}
We see that $ \left\{ \V_{v}  > x\right\} \subset
\Event{5}\cap\Event{6}\cap\Event{7}$ hence
\begin{equation} \label{tops1}
\Prob{\V_v  > x} \leq \Prob{\Event{5}} + \Prob{\Event{6}} + \Prob{\Event{7}}.
\end{equation}
We deal with each of  terms of the r.h.s. of (\ref{tops1}) separately. Let us choose $C
>\Cste{14}$ if $\S$ has
positive jumps and set $C=1$ otherwise. We now assume that $v$ and $a^{-1}(x)/v$ are very large.
According to (\ref{gd1}) and using the regular variation of
$a^{-1}(\cdot)$:
\begin{equation}
\label{tops2}\Prob{\Event{5}}\leq v\Prob{\V_1 >
(1-\varepsilon)x}\leq\frac{C}{(1-\varepsilon)^\alpha}\frac{v}{a^{-1}(x)}.
\end{equation}
We now deal with $\Prob{\Event{6}}$. Let $\eta>0$. Lemma \ref{regvar} gives for
all $v$ and $a^{-1}(x)/v$ large enough:
\begin{eqnarray*}
\frac{va^{-1}(x)}{(a^{-1}(z))^2} &=& \frac{a^{-1}\Pare{a(v)\frac{x}{a(v)}}}{a^{-1}(a(v))}
\Pare{\frac{a^{-1}(a(v))}{a^{-1}\Big(a(v)\Big(\frac{x}{a(v)}\Big)^\delta\Big)}}^2\\
&\leq & \Pare{\frac{x}{a(v)}}^{\alpha + \eta}\Pare{\frac{a(v)}{x}}^{2\delta(\alpha - \eta)}.
\end{eqnarray*}
Since $\delta > 1/2$, we can assume $\eta$ small enough such that
$2\delta(\alpha - \eta) - (\alpha + \eta) > \eta$ hence
\begin{equation}\label{gd3}
\frac{va^{-1}(x)}{(a^{-1}(z))^2} \leq \Pare{\frac{a(v)}{x}}^\eta,
\end{equation}
therefore, using (\ref{gd1}) then (\ref{gd3}):
\begin{equation}
\label{tops3}\Prob{\Event{6}} \leq v^2\Prob{\V_1 > z}^2\leq
C\frac{v^2}{(a^{-1}(z))^2}\leq
C\frac{v}{a^{-1}(x)}\Pare{\frac{a(v)}{x}}^\eta.
\end{equation}
Turning our attention to $\Prob{\Event{7}}$, we deduce from Tchebychev's inequality:
\begin{equation}
\label{tops4}\Prob{\Event{7}} \leq \frac{1}{\varepsilon^2
x^2}\Esp{(\zeta_{1,z} + \ldots + \zeta_{[v],z})^2} \leq
\frac{v}{\varepsilon^2 x^2}\Esp{\zeta_{1,z}^2} +
\frac{v^2}{\varepsilon^2 x^2}\Esp{\zeta_{1,z}}^2 .
\end{equation}
Let $f(z) = \Esp{(\zeta_{1,z})^2} = \int_{-z}^{z}y^2\Prob{\V_1\in
dy}$. This function is non decreasing and non zero for $z$ large
enough. It is also known from the characterization of the domain of
attraction (c.f. (8.14) of \cite{Felle1} p304) that the norming
constants $(a_n)$ are such that $nf(a_n)/a_n^2 \to \Cste{17} >0$,
hence $f(z) \sim \Cste{17}z^2/a^{-1}(z)$ as $z$ goes to infinity ($f$ is
regularly varying with index $2-\alpha$), therefore for $v$ and $a^{-1}(x)/v$
large enough:
\begin{equation}
\label{tops5}\frac{v}{\varepsilon^2 x^2}\Esp{(\zeta_{1,z})^2} = \frac{v
f(z)}{\varepsilon x^2} \leq
\Cste{18,\varepsilon}\frac{v}{a^{-1}(x)}\frac{f(z)}{f(x)} \leq
\Cste{18,\varepsilon}\frac{v}{a^{-1}(x)}.
\end{equation}
We can sharpen this estimate when $\alpha<2$. Indeed, in this case,
$f$ is regularly varying with index $2-\alpha >0$ thus using
Lemma \ref{regvar} and setting $\eta' =
(1-\delta)(2-\alpha)/2$:
\begin{equation*}
\frac{f(z)}{f(x)} \leq \Pare{\frac{z}{x}}^{(2-\alpha)/2} =
\Pare{\frac{a(v)(x/a(v))^\delta}{x}}^{(2-\alpha)/2} = \Pare{\frac{a(v)}{x}}^{\eta'}.
\end{equation*}
When $\alpha <2$, we therefore have:
\begin{equation}
\label{tops6}\frac{v}{\varepsilon^2 x^2}\Esp{(\zeta_{1,z}}\leq
\Cste{18,\varepsilon}\frac{v}{a^{-1}(x)}\Pare{\frac{a(v)}{x}}^{\eta'}.
\end{equation}
Let $g(z) = \Esp{\zeta_{1,z}} = \int_{-z}^{z}y\Prob{\V_1 \in dy}$.
Since $\V_1$ is in the domain of attraction of a stable law, it is
known that the centering constants $c(n)$ such that $\V_n/a(n) -
c(n)$ converge to a stable law may be chosen to be $c(n) = n
g(a(n))/a(n)$ (see \cite{Felle1} p305) but the assumption \ref{hyp1} of
this paper states that the norming constants $c(n)$ may also be
chosen to be $0$. This implies in particular that the sequence $n
g(a(n))/a(n)$ is bounded so we deduce that there 
exists $\Cste{19} > 0$ such that:
\begin{equation*}
|g(z)| \leq \Cste{19}\frac{z}{a^{-1}(z)} \hbox{ for all $z\geq 1$.}
\end{equation*}
Using this inequality, we get for $v$ and $a^{-1}(x)/v$ large enough:
\begin{eqnarray}
\nonumber\frac{v^2}{\varepsilon^2 x^2}\Esp{\zeta_{1,z}}^2 &\leq &\Cste{20,\varepsilon}\frac{v^2 z^2}{x^2 (a^{-1}(z))^2}\\
\nonumber& = & \Cste{20,\varepsilon}\frac{v}{a^{-1}(x)}\frac{v a^{-1}(x)}{(a^{-1}(z))^2}\Pare{\frac{z}{x}}^2\\
\nonumber& \leq &  \Cste{20,\varepsilon}\frac{v}{a^{-1}(x)}\frac{v a^{-1}(x)}{(a^{-1}(z))^2}\\
\label{tops7}& \leq &
\Cste{20,\varepsilon}\frac{v}{a^{-1}(x)}\Pare{\frac{a(v)}{x}}^\eta,
\end{eqnarray}
where we used (\ref{gd3}) for the last inequality. Putting the
pieces together,
(\ref{tops1})-(\ref{tops2})-(\ref{tops3})-(\ref{tops4})-(\ref{tops5})
and (\ref{tops7}) yield (\ref{regd1}). Moreover, when $\S$ has
positive jumps, we have $\alpha <2$, hence we can use (\ref{tops6})
instead of (\ref{tops5}) and we deduce that:
\begin{equation*}
\limsup_{\petit{
\begin{array}{c}
v\to\infty \\
\frac{a^{-1}(x)}{v}\to\infty
\end{array}%
}}\frac{a^{-1}(x)\Prob{\V_v > x}}{v} \leq \Cste{14}.
\end{equation*}
It remain to prove that the lower bound holds. Assume that $\S$
has positive jumps and notice that the event $\{\V_v > x\}$ contains
\begin{equation*}
\bigcap_{k=0}^{[v]-1}
\left\{\V^*_k\leq \varepsilon x\hbox{ , }
\V_{k+1}-\V_k > (1+2\varepsilon)x\hbox{ , }\Pare{\theta_{k+1}\V}^*_{[v]-k-1}\leq \varepsilon x\right\}.
\end{equation*}
Moreover the events of the last formulaare disjoints.
The independence and the stationarity of the increments of the random walk $\V$ yield
\begin{eqnarray*}
\Prob{\V_v > x}&\geq & \sum_{k=0}^{[v]-1} \P\Big(\V^*_k \leq \varepsilon x\Big)
\P\Big(\V_1 > (1+2\varepsilon)x\big)\P\Big(\V^*_{[v]-k-1} \leq \varepsilon x\Big)\\
&\geq& [v]\P\Big(\V^*_{v}\leq \varepsilon x\Big)^2\P\Big(\V_{1}>(1+2\varepsilon) x\Big).
\end{eqnarray*}
From (\ref{gd1}) and the regular variation of $a^{-1}(\cdot)$ we see that
$$[v]\Prob{\V_{1}>(1+2\varepsilon) x} \sim
\frac{\Cste{14}v}{{(1+2\varepsilon)}^\alpha a^{-1}(x)}$$
as $v$ and
$a^{-1}(x)/v$ both go to infinity. We also know from the results of
section \ref{subsectionD} that $\V^*_v / a(v)$ converges in
law towards $\S^*_1$ therefore:
\begin{equation*}
\lim_{\petit{
\begin{array}{c}
v\to\infty \\
\frac{a^{-1}(x)}{v}\to\infty
\end{array}%
}} \P\Big(\V^*_{v}\leq \varepsilon x\Big) = \lim_{\petit{
\begin{array}{c}
v\to\infty \\
\frac{a^{-1}(x)}{v}\to\infty
\end{array}%
}} \Prob{\frac{\V^*_{v}}{a(v)}\leq \varepsilon \frac{x}{a(v)}} = 1.
\end{equation*}
We conclude that
\begin{equation*}
\liminf_{\petit{
\begin{array}{c}
v\to\infty \\
\frac{a^{-1}(x)}{v}\to\infty
\end{array}%
}}\frac{a^{-1}(x)\Prob{\V_v > x}}{v} \geq \frac{\Cste{14}}{(1+2\varepsilon)^\alpha}.
\end{equation*}
\end{proof}

\begin{corr} \label{lemmefluc1}
By possibly extending the value of $\Cste{16}$, the equation
(\ref{regd1}) also holds with $\overline{\V}_v$, $-\underline{\V}_v$,
$\V^\#_v$ and $\V^*_v$ in place of $\V_v$.
\end{corr}
\begin{proof}
The results for $\overline{\V}_{v}$ and $-\underline{\V}_v$ are
straightforward  using Lemma \ref{lemmreflex}. As for $\V^*$ and
$\V^\#$, simply notice that $\{\V^\#_v \geq 2x \} \subset \{\V^*_v
\geq x\} \subset \{\overline{\V}_v \geq x\}\cup\{-\underline{\V}_v
\geq x\}$.
\end{proof}

\begin{corr} \label{cormomentsup}For any $0 < \delta < \alpha$:
\begin{equation*}
\lim_{v\rightarrow\infty}\Esp{\Pare{\frac{\overline{\V}_v}{a(v)}}^\delta}
= \Esp{\Pare{\overline{\S}_1}^\delta} \hbox{ and  }
\lim_{v\rightarrow\infty}\Esp{\left|\frac{\underline{\V}_v}{a(v)}\right|^\delta}
= \Esp{\Pare{-\underline{\S}_1}^\delta}.
\end{equation*}
\end{corr}
\begin{proof}
It follows from the last corollary and the regular variation of
$a^{-1}(\cdot)$ with index $\alpha$ that for any $0< \delta <
\alpha$:
\begin{equation*}
\sup_{v\geq
1}\Esp{\Pare{\frac{\overline{\V}_v}{a(v)}}^{\delta}}<\infty,
\end{equation*}
hence the family $\Pare{(\overline{\V}_v/a(v))^\delta , v\geq 1}$ is
uniformly integrable for all $0<\delta<\alpha$. We also know that
$\overline{\V}_v/a(v)$ converges in law toward $\overline{\S}_1$ as
$v$ goes to infinity. These two facts combined together yield the
first assertion. The proof of the second part of the corollary is
similar.
\end{proof}

\begin{prop} \label{lemmefluc3}
For all  $0 < \delta < q$ (recall that $q$ is the negativity
parameter of $\S$) there exists $\Cste{21,\delta}$ such that, for all
$v,x\geq 1$:
\begin{equation*}
\Prob{-\underline{\V}_{v} \leq x} \leq \Cste{21,\delta}
\Pare{\frac{a^{-1}(x)}{v}}^{\delta}.
\end{equation*}
We have a similar result for $\Prob{\overline{\V}_v \leq x}$ when
changing the condition $\delta < q$ by $\delta < p$.
\end{prop}

\begin{proof} We only prove the result for $\underline{\V}_v$.
By possibly extending the value of $\Cste{21,\delta}$, it suffice to prove
the inequality for $x$ and $v/a^{-1}(x)$ large
enough. Let us choose $\delta'$ such that $ \delta < \delta' < q <
1$ and notice that for any $y>0$:
\begin{eqnarray*}
\left\{ -\underline{\V}_{v} \leq x \right\} &\subset& \Lambda(x,y) \cup
\Pare{\left\{ -\underline{\V}_{v} \leq x\right\}\cap \Lambda(x,y)^c}\\
&\subset& \Lambda(x,y) \cup \left\{ \V^{\#}_{v} \leq x+y\right\},
\end{eqnarray*}
thus
\begin{equation} \label{inter20}
\Prob{-\underline{\V}_{v} \leq x} \leq \Prob{\Lambda\Pare{x,y}} +
\Prob{\V^{\#}_v \leq x+y}.
\end{equation}
On the one hand, for $x$ and $y/x$ large enough, using Proposition
\ref{tempssortie} and Lemma \ref{regvar}
\begin{eqnarray}
\nonumber\Prob{\Lambda\Pare{x,y}} &\leq& \Cste{22}\frac{b^{-1}\Pare{a^{-1}(x)}}{b^{-1}\Pare{a^{-1}(x+y)}} \\
\label{klmf1}&\leq&
\Cste{23,\delta'}\Pare{\frac{a^{-1}(x)}{a^{-1}(x+y)}}^{\delta'} .
\end{eqnarray}
On the other hand, for $x+y$ and  $v/a^{-1}(x+y)$ large enough,
using Proposition \ref{estimonte}:
\begin{eqnarray}
\label{klmf2}\Prob{\V^{\#}_v \leq x+y} &\leq &
\exp\Pare{-\frac{K^\#}{2}\frac{v}{a^{-1}(x+y)}}.
\end{eqnarray}
let us choose $y = a\Pare{\frac{K^\#
v}{2\log\Pare{v/a^{-1}(x)}}}-x$. It is easy to check that
(\ref{klmf1}) and (\ref{klmf2}) hold whenever $x$ and $v/a^{-1}(x)$
are large enough thus from (\ref{inter20}):
\begin{eqnarray*}
\Prob{-\underline{\V}_{v} \leq x} &\leq&
\Cste{23,\delta'}\Pare{\frac{2}{K^\#}}^{\delta'}\Pare{\frac{a^{-1}(x)}{v}
\Pare{\log\frac{v}{a^{-1}(x)}}}^{\delta '} + \frac{a^{-1}(x)}{v} \\
&\leq& \Cste{24,\delta'}\Pare{\frac{a^{-1}(x)}{v}}^{\delta} .
\end{eqnarray*}
\end{proof}

\section{Behavior of $X$}
In this section, we now study the diffusion $X$ in the random potential
$\V$. We will see that the behavior of this process depends strongly
on the environment. In order to do so, we will adapt the ideas of Hu and
Shi to our setting, in particular, we will show that the two Lemmas
$4.1$ and $4.2$ of \cite{HuShi1} still hold with a slight modification.

Recall  the well known  representation of $X$ (c.f.
\cite{Brox1,HuShi1,Ito1}) which states that we can construct $X$ from
a Brownian motion through a (random) change of scale and a (random)
change of time hence we will assume that $X$ has the form:
\begin{equation} \label{representationX}
X_t = \Achange^{-1}\Pare{B_{\Tchange^{-1}(t)}}
\end{equation}
where $B$ is a standard Brownian motion independent of $\V$ and
where $\Achange^{-1}$ and  $\Tchange^{-1}$ are the respective
inverses of
\begin{eqnarray*}
\Achange\Pare{x} &=& \int_{0}^{x}e^{\V_y}dy  \qquad \hbox{ for $x\in\R,$}\\
\Tchange\Pare{t} &=&
\int_{0}^{t}e^{-2\V_{\Achange^{-1}\Pare{B_s}}}ds \qquad \hbox{ for
$t\geq 0.$}
\end{eqnarray*}
Note that our  assumption on $\V$ implies with probability $1$ that
$\Achange$ is an increasing homeomorphism on $\R$ and that $\Tchange$
is an increasing homeomorphism on $\R_+$, thus $\Achange^{-1}$ and
$\Tchange^{-1}$ are well defined. Let $v>0$ and recall the definition
of $\sigma$ given in section \ref{subsectionD}. Using
(\ref{representationX}) we have:
\begin{equation*}
\sigma_{X}(v) = \Tchange\Pare{\sigma_{B}\Pare{\Achange(v)}}.
\end{equation*}
Let $(L(t,x) , t\geq 0 , x\in\R)$ stand for the bicontinuous version
of the local time process of $B$. The last equality may be
rewritten:
\begin{eqnarray*}
\sigma_{X}(v) &=& \int_{0}^{\sigma_{B}(\Achange(v))}e^{-2\V_{\Achange^{-1}(B_s)}}ds \\
&=& \int_{-\infty}^{\Achange(v)} e^{-2\V_{\Achange^{-1}(x)}}L(\sigma_{B}(\Achange(v)),x)dx \\
&=&  \int_{-\infty}^{v}
e^{-\V_{y}}L(\sigma_{B}(\Achange(v)),\Achange(y))dy
\end{eqnarray*}
where we have used the change of variable $x = \Achange(y)$. Let us now define $I_1$ and $I_2$:
\begin{eqnarray}
\label{expI1}I_1 (v) &=& \int_{0}^{v}e^{-\V_{y}}L(\sigma_{B}(\Achange(v)),\Achange(y))dy, \\
\label{expI2}I_2 (v) &=& \int_{0}^{\infty}e^{-\V_{-y}}L(\sigma_{B}(\Achange(v)),\Achange(-y))dy.
\end{eqnarray}
Using the definition of $\sigma_{X}$, we get
\begin{equation} \label{egaliteset}
\left\{\overline{X}_t \geq v\right\} = \left\{ I_1(v) + I_2(v) \leq t\right\}.
\end{equation}
The next two propositions show the connection between $\V$ and $X$.
These estimates will enable us to reduce the study of the limiting
behavior of $X$ to the study of some functionals of the potential
$\V$. The streamline of the proofs is the same as that of Lemmas
$4.1$ and $4.2$ of \cite{HuShi1} and one should refer to the proof
of these two lemmas for further details.
\begin{prop}\label{propI1} there exists $\Cste{25}$ such that for all $v$ large enough
\begin{equation*}
\V^{\#}_{v-\frac{1}{2}} - \Pare{\log v}^{4} \leq \log I_1 (v) \leq
\V^{\#}_{v} + \Pare{\log v}^{4} \hbox{ on $\Event{8}(v),$}
\end{equation*}
where $\Event{8}(v)$ is a measurable set such that
\begin{equation*}
\Prob{\Event{8}(v)^{c}}\leq \Cste{25}e^{-(\log v)^{2}}.
\end{equation*}
\end{prop}

\begin{prop}\label{propI2} there exists $\Cste{26}$ such that for all $v$ large enough
\begin{eqnarray*}
\log I_2(v) &\leq& \widetilde{U}_{\V}\Pare{\overline{\V}_v + (\log v)^{4}} \hbox{ on $\Event{9}(v),$}\\
\log I_2(v) &\geq&
\widetilde{U}_{\V}\Pare{\overline{\V}_{v-\frac{1}{2}} - (\log
v)^{4}} \hbox{ on $\Event{9}(v) \cap \left\{
\overline{\V}_{v-\frac{1}{2}} > (\log v)^{4}\right\},$}
\end{eqnarray*}
where $\widetilde{U}$ was defined in section \ref{subsectionD} and where
$\Event{9}(v)$ is a measurable set such that
\begin{equation*}
\Prob{\Event{9}(v)^{c}}\leq \Cste{26}e^{-(\log v)^{2}}.
\end{equation*}
\end{prop}

\begin{proof}[Proof of Proposition \ref{propI1}.] For $v>0$, let $\Rbes^2$ be defined as:
\begin{eqnarray*}
\Rbes^2 (t) = \frac{L\Pare{\sigma_{B}(\Achange(v)) , \Achange(v)
-t\Achange(v)}}{\Achange(v)} \qquad \hbox{ for $0\leq t\leq 1$}.
\end{eqnarray*}
Let $\Rbes$ be the positive root of $\Rbes^2$. Just as in
\cite{HuShi1} , p1498, we see, using Ray-Knight Theorem and the
scaling property of the Brownian motion that for any fixed $v$ the
process $(\Rbes(t),0\leq t\leq 1)$ has the law of a two dimensional
Bessel process starting from $0$. Moreover, $\Rbes$ is independent
of $\V$. We can now rewrite (\ref{expI1}) as
\begin{equation*}
I_1(v) = \Achange(v) \int_{0}^{v} e^{-\V_s}\Rbes^2\left(
\frac{\Achange(v) -\Achange(s)}{\Achange(v)}\right)ds.
\end{equation*}
Let us define
$$\Event{10} = \left\{ \sup_{0<t\leq 1}\frac{\Rbes(t)}{\sqrt{t\log(8/t)}}\leq \sqrt{v}\right\}.$$
Using Lemma $6.1$ p1497 of \cite{HuShi1}, we get
$\Prob{\Event{10}^c}\leq \Cste{27} e^{-v/2}$. On $\Event{10}$, we have
$$I_1(v) \leq
v\int_{0}^{v}e^{-\V_s}\left(\Achange(v)-\Achange(s)\right)
\log\left(\frac{8\Achange(v)}{\Achange(v)-\Achange(s)}\right)ds,$$
and for all $s\leq v$
\begin{equation*}e^{-\V_s}\left(\Achange(v)-\Achange(s)\right) =
\int_{s}^{v}e^{\V_y-\V_s}dy \leq v e^{\V^{\#}_v}.\end{equation*}
This implies:
\begin{equation} \label{oldeq1}
I_1(v) \leq v^2
e^{\V^\#_v}\int_{0}^{v}\log\left(\frac{8\Achange(v)}{\Achange(v)-\Achange(s)}\right)ds.
\end{equation}
We also have
$$\Achange(v) = \int_{0}^{v} e^{\V_s}ds \leq v e^{\overline{\V}_v}\quad\hbox{ and }\quad \Achange(v)
- \Achange(s) = \int_{s}^{v}e^{\V_y}dy \geq (v-s)e^{\underline{\V}_v},$$
thus
\begin{eqnarray*}
\int_{0}^{v}\log\left(\frac{8\Achange(v)}{\Achange(v)-\Achange(s)}\right)ds
&\leq& v\left(\overline{\V}_v -
\underline{\V}_v\right)+\int_{0}^{v}\log\left(\frac{8v}{v-s}\right)ds \\
&\leq &v\left( \overline{\V}_v - \underline{\V}_v + 1 +\log(8)\right).
\end{eqnarray*}
Combining this with (\ref{oldeq1}) yields $\log(I_1(v)) \leq \V^\#_v
+ \log\left(\overline{\V}_v - \underline{\V}_v\right) + 4\log(v)$
for all $v$ large enough. We now define $\Event{11} =
\left\{ \log\left(\overline{\V}_v - \underline{\V}_v\right)\leq
\log^3 (v)\right\}$. On  $\Event{10}\cap \Event{11}$, for all $v$ large
enough, we get the upper bound :
$$\log(I_1(v)) \leq \V^\#_v + \log^4 (v).$$
Notice that $\{\overline{\V}_v - \underline{\V}_v > a\}\subset\{
\V^{*}_v > a/2\}$ thus using Corollary \ref{lemmefluc1} and the
regular variation of $a^{-1}(\cdot)$, it is easily checked that
$\Prob{\Event{11}^c} \leq \exp(-\log^2(v))$ for any $v$ large enough.
We now prove the existence of the lower bound. For the sake of clarity,
we will use the notation $l = \log(v)$ and $\delta = \exp(-l^2)$.
For $v>1/2$, there exist two integers $0\leq k^- \leq k^+ \leq
v-\frac{1}{2}$ such  that $\V^{\#}_{v-\frac{1}{2}} = \V_{k^+} -
\V_{k^-}$. Let us define the sets:
\begin{eqnarray*}
\Event{12} &=& \left\{ \inf_{k^-\leq s\leq k^- +
\frac{1}{2}}\Rbes\left(\frac{\Achange(v)-\Achange(s)}{\Achange(v)}\right)
>
\delta \sqrt{\frac{\Achange(v)-\Achange(k^-)}{\Achange(v)}}\right\},  \\
\Event{13} &=& \left\{ \V^{\#}_{v-\frac{1}{2}}\geq 3 l^2\right\}.
\end{eqnarray*}
Using again Lemma $6.1$ p1497 of \cite{HuShi1} combined
with the independence of $\Rbes$ and $\V$:
\begin{equation}\label{oldeq2}
\Prob{(\Event{12}\cap \Event{13})^c} \leq \Prob{\Event{13}^c} +
2\delta +
2\Esp{e^{-\frac{\delta^2}{2}J(v)}\Indic_{\Event{13}}},\end{equation}
where $J$ is given by:
$$J(v) = \frac{\Achange(v) - \Achange(k^-)}{\Achange\Pare{k^- +
\frac{1}{2}}-\Achange(k^-)}.$$
On the one hand:
$$\Achange(v) - \Achange(k^-) = \int_{k^-}^{v}e^{\V_{s}}ds
\geq \int_{k^+}^{k^+ + \frac{1}{2}}e^{\V_{s}}ds =
\frac{1}{2}e^{\V_{k^+}}.$$
On the other hand, since $k^-$ is a
integer and $\V$ is flat on $[k^-,k^-+1)$ we have:
$$\Achange\Pare{k^- + \frac{1}{2}} - \Achange(k^-) = \int_{k^-}^{k^- + \frac{1}{2}}e^{\V_{s}}ds
= \frac{1}{2}e^{\V_{k^-}}.$$ This implies $J(v) \geq
\exp(\V^\#_{v-1/2})$. Using this inequality combined with (\ref{oldeq2}), we
get:
$$ \Prob{(\Event{12}\cap \Event{13})^c} \leq \Prob{\Event{13}^c}
+ 2\delta + 2\exp(-\delta^2 \exp(3 l^2)/2),$$ hence for any $v$ large
enough, we have $\Prob{(\Event{12}\cap \Event{13})^c} \leq
\Prob{\Event{13}^c} + 3\exp(-l^2).$ Using Proposition
\ref{estimonte}, it is easily seen that $\Prob{\Event{13}^c} \leq
e^{-l^2}$ for all large enough $v$'s. Let us finally set $\Event{8} =
\Event{10} \cap \Event{11} \cap \Event{12} \cap \Event{13}$. We have
proved that there exists $\Cste{25}>0$ such that $\Prob{\Event{8}^c}
\leq \Cste{25} \exp(-l^2)$. Notice that:
\begin{eqnarray*}
I_1(v) &=& \Achange(v) \int_{0}^{v} e^{-\V_s}\Rbes^2\left(
\frac{\Achange(v) -\Achange(s)}{\Achange(v)}\right)ds \\
&\geq&  \Achange(v)
e^{-\V_{k^-}} \int_{k^-}^{k^- + \frac{1}{4}} \Rbes^2\left(
\frac{\Achange(v) -\Achange(s)}{\Achange(v)}\right)ds,
\end{eqnarray*}
therefore on $\Event{8}$:
\begin{equation*}
I_1(v) \geq  \delta^2 e^{-\V_{k^-}} \int_{k^-}^{k^- + \frac{1}{4}}\left(\Achange(v)
-\Achange(s)\right) ds,
\end{equation*}
 but for all $s$ such that  $k^- \leq s \leq k^- + \frac{1}{4}$ we also have
\begin{equation*}
\Achange(v) - \Achange(s)\geq \Achange(v) - \Achange\left(k^- +
\frac{1}{4}\right) = \int_{k^- + \frac{1}{4}}^{v}e^{\V_y}dy\geq
\int_{k^+ +\frac{1}{4}}^{k^+
+\frac{1}{2}}e^{\V_y}dy=\frac{1}{4}e^{\V_{k^+}},
\end{equation*}
hence
\begin{equation*}
\int_{k^-}^{k^- + \frac{1}{4}}\left(\Achange(v)
-\Achange(s)\right) ds \geq \frac{1}{16}e^{\V_{k^+}}.
\end{equation*}
We finally get on $\Event{8}$:
\begin{equation*}
I_1(v) \geq \frac{\delta^2}{16}e^{\V^\#_{v -\frac{1}{2}}}.
\end{equation*}
We conclude the proof of the proposition by taking the logarithm.
\end{proof}

\begin{proof}[Proof of Proposition \ref{propI2}.]
For $v> 0$, we define the process $\Zbes$ by
\begin{equation*}
\Zbes(t) = \frac{L\Pare{\sigma_{B}\Pare{\Achange(v)} ,
-t\Achange(v)}}{\Achange(v)} \quad \hbox{ for $t\geq 0.$}
\end{equation*}
Using Ray-Knight Theorem and the scaling property of the Brownian
motion, we see that for any fixed $v$ the process $\Zbes$ has the
law of a squared Bessel process such that $\Zbes(0)$ has an
exponential distribution with mean $2$. Moreover, $\Zbes$ is
independent of $\V$. We can now rewrite (\ref{expI2}):
\begin{equation*}
I_2 (v) = \Achange(v)\int_{0}^{\infty}e^{-\V_{-s}}\Zbes\left(\frac{-\Achange(-s)}{\Achange(v)}\right)ds.
\end{equation*}
We know that $0$ is an absorbing state for $\Zbes$. Let $\zeta =
\inf\Pare{s\geq 0 \hbox{ , } \Zbes_s = 0}$ be the absorption time of
$\Zbes$ and let us also define
\begin{equation*}
\zeta(v) = \inf\left(s\geq 0\hbox{ ,
}\Zbes\left(\frac{-\Achange(-s)}{\Achange(v)}\right)=0 \right).
\end{equation*}
We can now write
\begin{equation*}
I_2 (v) = \Achange(v)\int_{0}^{\zeta(v)}e^{-\V_{-s}}\Zbes\left(\frac{-\Achange(-s)}{\Achange(v)}\right)ds.
\end{equation*}
We keep the notation $l = \log(v)$, note that $\Achange(v) =
\int_{0}^{v}e^{\V_s}ds \leq \exp(\overline{\V}_{v} + l)$, therefore
\begin{eqnarray*}
I_2 (v) &\leq& e^{\overline{\V}_{v} +
l}\zeta(v)\sup_{0\leq s\leq
\zeta(v)}\left(e^{-\V_{-s}}\right)\sup_{s\geq 0}\Zbes(s) \\
&\leq& \zeta(v)\sup_{s\geq 0}\Zbes(s)e^{l + \overline{\V}(v) -
\underline{\V}(-\zeta(v))}.
\end{eqnarray*}
Let us define $\Event{14} = \{ \sup_{s\geq 0} \Zbes(s) \leq
\exp(l^2)\}$, using Lemma $7.1$, p1501 of \cite{HuShi1}, we find:
$\P(\Event{14}^c)\leq 4 \exp(-l^2)$, thus on $\Event{14}$, we have:
\begin{equation}\label{oldeq3} I_2(v) \leq
\zeta(v)e^{2l^2+\overline{\V}(v) - \underline{\V}(-\zeta(v))}.
\end{equation}
Let $\Event{15} = \left\{ \zeta(v) \leq
\widetilde{\sigma}_{\V}\left(\overline{\V}_v + l^4\right)
+\frac{1}{2}\right\}$ and notice that for all $a\geq 0$:
$$\left\{ \zeta(v)>a\right\} =\left\{\frac{-\Achange(-a)}{\Achange(v)}<\zeta\right\}.$$ Therefore
\begin{equation*}
\Prob{\Event{15}^c}  = \Prob{ \frac{ -\Achange\left(
-\widetilde{\sigma}_{\V}\Pare{\overline{\V}_v + l^4}-\frac{1}{2}
\right)} {\Achange(v)}<\zeta },
\end{equation*}
but
\begin{equation*}
-\Achange\Pare{-\widetilde{\sigma}_{\V}(\overline{\V}_v +
l^4)-\frac{1}{2}} \geq
\int_{-\widetilde{\sigma}_{\V}(\overline{\V}_v +
l^4)-\frac{1}{2}}^{-\widetilde{\sigma}_{\V}(\overline{\V}_v +
l^4)}e^{\V_s}ds\geq \frac{1}{2}e^{\overline{\V}_v + l^4},
\end{equation*}
and we have already seen that $\Achange_{v}
\leq \exp(\overline{\V}_{v} + l)$, combining this two
inequalities yields for all large enough $v$'s:
\begin{equation*}
\frac{-\Achange\Pare{-\widetilde{\sigma}_{\V}(\overline{\V}_v +
l^4)-\frac{1}{2}}}{\Achange(v)}\geq e^{l^3},
\end{equation*}
hence
\begin{equation*}
\Prob{\Event{15}^c} \leq \Prob{\zeta > e^{l^3}}\leq e^{-l^3},
\end{equation*}
where we have used Lemma $7.1$ p1501 of \cite{HuShi1} for the last
inequality. On $\Event{14}\cap \Event{15}$, for $v$ large enough, we
deduce from (\ref{oldeq3}) the inequality:
\begin{equation*}
I_2(v) \leq \zeta(v)e^{2l^2+\overline{\V}(v) -
\underline{\V}(-\widetilde{\sigma}_{\V}\left(\overline{\V}_v +
l^4\right) +\frac{1}{2})}.
\end{equation*}
But $\overline{\V}(v) -
\underline{\V}(-\widetilde{\sigma}_{\V}\left(\overline{\V}_v +
l^4\right) +\frac{1}{2}) = \widetilde{U}_{\V}(\overline{\V}_v +l^4)
-l^4$ (recall that $\V$ is flat on $(-n-1,-n] \hbox{ , } n \in N$).
Therefore on  $\Event{14} \cap \Event{15}$:
\begin{equation*}
I_2(v) \leq \zeta(v) e^{-l^3
+\widetilde{U}_{\V}\left(\overline{\V}_v + l^4\right)}.
\end{equation*}
Let $\Event{16} = \{ \widetilde{\sigma}_{\V}(\overline{\V}_{v} +
l^4)+ \frac{1}{2}\leq \exp(l^3)\}$. On $\Event{17} = \Event{14} \cap
\Event{15} \cap \Event{16}$ we have $\zeta(v) \leq \exp(l^3)$, hence
on $\Event{17}$, for all large enough $v$'s:
\begin{equation*}
\log\left(I_2 (v)\right) \leq \widetilde{U}_{\V}\left(
\overline{\V}_{v} + \log^4 v\right),
\end{equation*}
this gives the upper bound on $\Event{17}$. Let us check that
$\Prob{\Event{16}^c} \leq \Cste{28} \exp(-l^2)$. We have $\P(\Event{16}^c)\leq
\Prob{\widetilde{\sigma}_{\V}(\overline{\V}_{v} + l^4) >
\exp(l^3)/2}$ thus
\begin{equation*}\label{al89}
\Prob{\Event{16}^c}\leq \Prob{\overline{\V}(-\frac{1}{2}e^{l^3})\leq
2\overline{\V}(v)}+\Prob{\overline{\V}(-\frac{1}{2}e^{l^3})\leq 2
l^4}.
\end{equation*}
We also have
\begin{equation*}
\Prob{\overline{\V}(-\frac{1}{2}e^{l^3})\leq2\overline{\V}(v)} \leq
\Prob{\overline{\V}(-\frac{1}{2}e^{l^3})\leq e^{l^{5/2}}} +
\Prob{\overline{\V}(v) > \frac{1}{2}e^{l^{5/2}}}.
\end{equation*}
Using Corollary \ref{lemmefluc1} and the regular variation of
$a^{-1}(\cdot)$, for all $v$ large enough:
\begin{equation*}
\Prob{\overline{\V}(v) > \frac{1}{2}e^{l^{5/2}}} \leq e^{-l^2}.
\end{equation*}
Recall that $(\V(x) , x\geq 0)$ and $(-\V(-x),x \geq 0)$ have the
same law thus Proposition \ref{lemmefluc3} implies:
\begin{equation*}
\Prob{\overline{\V}(-\frac{1}{2}e^{l^3})\leq 2 l^4} \leq
\Prob{\overline{\V}(-\frac{1}{2}e^{l^3})\leq e^{l^{5/2}}} \leq
e^{-l^2}.
\end{equation*}
These inequalities give $\Prob{\Event{16}^c} \leq 3e^{-l^2}$ hence
$\Prob{\Event{17}^{c}} \leq 8 e^{-l^{2}}$. We now prove the lower bound.
Notice that
\begin{equation}
\label{oldeq12} \Achange(v) \geq
\int_{\sigma_{\V}(\overline{\V}(v-\frac{1}{2}))}^{\sigma_{\V}(\overline{\V}(v-\frac{1}{2}))+\frac{1}{2}}
e^{\V_s}ds= \frac{1}{2}e^{\overline{\V}(v-\frac{1}{2})},
\end{equation}
and for all $x\leq
\widetilde{\sigma}_{\V}(\overline{\V}_{v-\frac{1}{2}}-l^4 ) \leq
\widetilde{\sigma}_{\V}(\overline{\V}_v)$:
\begin{equation}
\label{oldeq12zzz}-\Achange(-x)  = \int_{-x}^{0}e^{\V(s)}ds \leq
e^{\overline{\V}(v-\frac{1}{2})-l^4}\widetilde{\sigma}_{\V}(\overline{\V}_{v}),
\end{equation}
therefore, for all  $x\leq
\widetilde{\sigma}_{\V}(\overline{\V}(v-\frac{1}{2})-l^4)$ we have
$-\Achange_{-x}/\Achange_{v} \leq
\exp(-l^4)\widetilde{\sigma}_{\V}(\overline{\V}_v)$. Let $\Event{18}
= \{ \widetilde{\sigma}_{\V}(\overline{\V}_{v})\leq \exp(l^3)\}$. As for
the estimate of $\Prob{\Event{16}^c}$, it is easily checked that for all
$v$ large enough, $\Prob{\Event{18}^c}\leq 3 \exp(-l^2)$. Moreover,
on the set $\Event{18}$, combining (\ref{oldeq12}) and (\ref{oldeq12zzz}),
we have $-\Achange(-x)/\Achange(v) \leq
e^{-\frac{1}{2}l^4}$ for all $0\leq x\leq
\widetilde{\sigma}_{\V}(\overline{\V}(v-\frac{1}{2})-l^4)$. Let us
now define
\begin{equation*}
\Event{19} = \left\{\inf_{0\leq s\leq e^{-\frac{1}{2}l^4}}\Zbes(s)
\geq e^{-l^2}\right\}.
\end{equation*}
Using Lemma $7.1$ p1501 of \cite{HuShi1}, we see that
$\P(\Event{19}^c)\leq 2 e^{-l^2}$. Recall
that:
\begin{eqnarray*}
I_2(v)&=&\Achange(v)\int_{0}^{\infty}e^{-\V_{-s}}\Zbes\left(\frac{-\Achange(-s)}{\Achange(v)}\right)ds \\
&\geq&
\Achange(v)\int_{0}^{\widetilde{\sigma}_{\V}(\overline{\V}(v-\frac{1}{2})-l^4)}
e^{-\V_{-s}}\Zbes\left(\frac{-\Achange(-s)}{\Achange(v)}\right)ds,
\end{eqnarray*}
therefore on $\Event{20} = \Event{18}\cap \Event{19}$:
\begin{equation*}
I_2(v) \geq
\widetilde{\sigma}_{\V}\Big(\overline{\V}\big(v-\frac{1}{2}\big)-l^4\Big)
\Achange(v)
e^{-\underline{\V}\Pare{-\widetilde{\sigma}_{\V}(\overline{\V}(v-\frac{1}{2})-l^4)}-l^2}.
\end{equation*}
Using again (\ref{oldeq12}) we find on $\Event{20}$:
\begin{eqnarray*}
I_2(v) &\geq& \frac{1}{2}
\widetilde{\sigma}_{\V}\Big(\overline{\V}\big(v-\frac{1}{2}\big)-l^4\Big)
e^{\overline{\V}(v-\frac{1}{2})-\underline{\V}\Pare{-\widetilde{\sigma}_{\V}(\overline{\V}(v-\frac{1}{2})-l^4)}-l^2}\\
&=& \frac{1}{2}
\widetilde{\sigma}_{\V}\Big(\overline{\V}\big(v-\frac{1}{2}\big)-l^4\Big)
e^{\widetilde{U}_{\V}\Pare{\overline{\V}(v-\frac{1}{2})-l^4}+l^4-l^2}.
\end{eqnarray*}
Notice that on  $\{ \overline{\V}\big(v-1/2\big) > l^4 \}$, we
have $\widetilde{\sigma}_{\V}(\overline{\V}(v-1/2)-l^4 )\geq 1$ (because $\V$ is identically $0$ on $(-1,0]$).
This implies that on $\Event{20}\cap \{
\overline{\V}(v-1/2) > l^4 \}$:
\begin{equation*}
I_2(v) \geq
e^{\widetilde{U}_{\V}\Pare{\overline{\V}\big(v-\frac{1}{2}\big) -
l^4}},
\end{equation*}
which yields the lower bound by taking the logarithm. Finally, let
$\Event{9} = \Event{20}\cap \Event{17}$, we have
\begin{equation*}
\Prob{\Event{9}^c} \leq \Prob{\Event{17}^c} + \Prob{\Event{20}^c}
\leq 13 e^{-(\log v)^2}
\end{equation*}
for all large enough $v$'s and the upper bound holds on $\Event{9}$ as
well as the lower bound on $\Event{9}\cap \{\overline{\V}\big(v-1/2\big)> l^4 \}$.
\end{proof}

\section{Proof of the main theorems}

\subsection{Proof of Theorem \ref{MainTheo1}}
We first state two lemmas before we give the proof of the theorem.
\begin{lemm} \label{lemme1theo1} For any $c_0 >0$, we have
\begin{eqnarray*}
\limsup_{t\rightarrow\infty}\frac{\log \Prob{\overline{X}_t \geq c_0
a^{-1}\Pare{\log t}\log\log\log t}}{\log\log\log t} \leq -c_0
K^\#,
\end{eqnarray*}
where $K^\#$ was defined in Proposition \ref{estimonte}.
\end{lemm}
\begin{proof} Let $v = c_0 a^{-1}\Pare{\log t}\log\log\log t$, using
(\ref{egaliteset}) and Proposition \ref{propI1} we get for all $t$
large enough:
\begin{eqnarray*}
\Prob{\overline{X}_t \geq v} &\leq& \Prob{I_1(v) \leq t}\\
&\leq & \Prob{\V^{\#}_{v-\frac{1}{2}} \leq \log t + (\log v)^{4}} +
\Cste{25}\exp\Pare{-(\log v)^{2}}.
\end{eqnarray*}
Using Proposition \ref{estimonte} , for any $\varepsilon>0$ and for
all $t$ large enough (depending on $\varepsilon$):
\begin{eqnarray*}
\Prob{\V^{\#}_{v-\frac{1}{2}} \leq \log t + (\log v)^{4}} & \leq &
\exp\Pare{-(K^\# -\varepsilon)\frac{v-1/2}{a^{-1}\Pare{\log t + (\log v)^4}}} \\
&\leq& \exp\Pare{-c_0 (K^\# -2\varepsilon)\log\log\log t}
\end{eqnarray*}
where we used the regular variation of $a^{-1}(\cdot)$ to check that
$a^{-1}(\log t + (\log v)^4)\sim a^{-1}(\log t)$. We therefore
obtain for all $t$ large enough:
\begin{eqnarray*}
\Prob{\overline{X}_t \geq v} &\leq&  \exp\Pare{-c_0 (K^\#
-2\varepsilon)\log\log\log t} + \exp\Pare{-(\log v)^2}\\
&\leq&2\exp\Pare{-c_0 (K^\# -2\varepsilon)\log\log\log t}.
\end{eqnarray*}
\end{proof}
\begin{lemm} \label{lemme2theo1}
For any $c_0 >0$ and for all $t$ large enough (depending on $c_0$) we have
\begin{eqnarray*}
\left\{\overline{X}_t \geq v\right\} &\supset& \left\{\V^{\#}_{v} \leq
\log t - \sqrt{\log t} \hbox{ , } \overline{\V}_v \leq \frac{\log
t}{5}\right\}\\
&&\cap\left\{\widetilde{U}_{\V}\Pare{\frac{\log
t}{4}}\leq \frac{\log t}{2}\right\}\cap \Event{21}(v)
\end{eqnarray*}
where $v = c_0 a^{-1}(\log t)\log\log\log t$ and where
$\Event{21}(v)$ is a measurable set such that:
\begin{equation*}
\Prob{\Event{21}^c(v)} \leq \Cste{29}e^{-(\log v)^{2}}.
\end{equation*}
\end{lemm}
\begin{proof}
Using (\ref{egaliteset}) combined with Proposition \ref{propI1} and
\ref{propI2} we get, for $t$ sufficiently large
\begin{eqnarray*}
\left\{\overline{X}_t \geq v \right\} &=& \left\{I_1 (v) + I_2 (v)
\leq t\right\}\\
&\supset &\left\{e^{\V^{\#}_v + (\log v)^{4}} +
e^{\widetilde{U}_{\V}\Pare{\overline{\V}_v + (\log v)^{4}}} \leq
t\right\}\cap \Event{21}(v)
\end{eqnarray*}
with $\Event{21}(v) = \Event{8}(v) \cap \Event{9}(v)$ thus
$\Prob{\Event{21}^{c}(v)} \leq \Cste{29}e^{-(\log v)^2}$. Notice
also that
\begin{equation*}
\left\{ \V^{\#}_v \leq \log t - \sqrt{\log t}\right\} \subset
\left\{ \V^{\#}_v + \log^{4} v \leq \log\frac{t}{2}\right\},
\end{equation*}
hence $\left\{\overline{X}_t \geq v \right\}$ contains
\begin{equation}\label{baba1}
\left\{ \V^{\#}_v \leq
\log t - \sqrt{\log t}\right\} \cap
\left\{\widetilde{U}_{\V}\Pare{\overline{\V}_v + (\log v)^{4}} \leq
\log\Pare{\frac{t}{2}}\right\}\cap \Event{21}(v).
\end{equation}
We also have $\left\{\overline{\V}_v \leq \frac{\log t}{5}\right\}
\subset \left\{\overline{\V}_v +(\log v)^{4} \leq \frac{\log
t}{4}\right\}$ therefore:
\begin{equation*}
\left\{\overline{\V}_v \leq \frac{\log t}{5} \hbox{ , }
\widetilde{U}_{\V}\Pare{\frac{\log t}{4}} \leq \frac{\log
t}{2}\right\} \subset \left\{
\widetilde{U}_{\V}\Pare{\overline{\V}_v + (\log v)^{4}} \leq
\frac{\log t}{2}\right\},
\end{equation*}
combining this with (\ref{baba1}) completes the proof.
\end{proof}

\begin{proof}[Proof of Theorem \ref{MainTheo1}]
As we already mentioned in the introduction, $X$ and $\overline{X}$
have the same  upper function so we only need to prove the theorem for
$\overline{X}$. Let us choose $K$ such that $K < K^\#$ and
$\varepsilon >0$. Define the sequence $t_i = \exp(\exp(\varepsilon
i))$. We also use the notation $f(x) = a^{-1}(\log x)\log\log\log
x$. Using regular variation of $a(\cdot)$ we easily check that
$f(t_i) /f(t_{i+1})$ converges to $\exp(-\alpha\varepsilon)$ thus,
for all $i$ large enough
\begin{equation*}
\Prob{\overline{X}_{t_{i+1}}\geq \frac{f(t_i)}{K}} \leq
\Prob{\overline{X}_{t_{i+1}}\geq
\frac{f(t_{i+1})}{e^{2\varepsilon}K}}.
\end{equation*}
Using Lemma \ref{lemme1theo1}:
\begin{equation*}
\limsup_{i\rightarrow\infty}\frac{1}{\log(\varepsilon(i+1))}\log
\Pare{\Prob{\overline{X}_{t_{i+1}}\geq \frac{f(t_i)}{K}}} \leq
-\frac{K^\#}{e^{2\varepsilon}K}.
\end{equation*}
Since $K < K^\#$, we can choose $\varepsilon$ small enough such that
$K^\# / (K\exp(2\varepsilon)) < 1$ and we deduce from the last
inequality that the sum $\sum \P(\overline{X}_{t_{i+1}}\geq
f(t_i)/K)$ converges. Using Borel-Cantelli Lemma, with probability
$1$, for all $i$ large enough $\overline{X}_{t_{i+1}}\leq f(t_i)/K$.
For $t\in [t_i , t_{i+1}]$, using monotonicity of $f$ and
$\overline{X}$:
\begin{equation*}
\overline{X}_t \leq \overline{X}_{t_{i+1}} \leq \frac{f(t_{i})}{K}
\leq \frac{f(t)}{K}.
\end{equation*}
This holds for all $K<K^\#$ hence we proved that almost surely:
$$\limsup_{t\rightarrow\infty}\frac{\overline{X}_t}{f(t)}\leq \frac{1}{K^\#}$$. We
now prove the lower bound. Choose $K > K^\#$ and change the sequence
$(t_i)$ for $t_i = \exp(\exp i)$. From Lemma \ref{lemme2theo1}, for
$i$ large enough:
\begin{equation*}
\left\{\overline{X}_{t_i}\geq \frac{f(t_i)}{K}\right\} \supset
\Event{21}(f(t_i)/K) \cap \Event{22}(i),
\end{equation*}
where $\Event{21}$ was defined in Lemma \ref{lemme2theo1} and where
$\Event{22}(i) = \Event{23}(i)\cap\Event{24}(i)\cap\Event{25}(i)$
with
\begin{eqnarray*}
\Event{23}(i) &=& \left\{ \widetilde{U}_{\V}\Pare{e^i/4}\leq e^i/2\right\},\\
\Event{24}(i) &=& \left\{ \V^{\#}_{f(t_i)/K} \leq e^i -
e^{i/2}\right\},\\
\Event{25}(i) &=& \left\{ \overline{\V}_{f(t_i)/K} \leq e^i
/5\right\}.
\end{eqnarray*}
Moreover, $\sum \P(\Event{21}^c(f(t_i)/K))<\infty$ so it only
remains to be proved that the events $\Event{22}(i)$ happen
infinitely often almost surely. It follows from results of section \ref{subsectionD} that
$\lim_{i\rightarrow\infty}\P(\Event{23}(i)) =
\P(\widetilde{U}_{\S}(1/4)\leq 1/2)$ and it is clear that this quantity is
not $0$. Since $\Event{24}(i)\cap\Event{25}(i)$ and $\Event{23}(i)$
are independent events $\P(\Event{22}(i)) \geq
\Cste{30}\P(\Event{24}(i)\cap\Event{25}(i))$ for all $i$ large
enough thus, we deduce from Proposition \ref{Encadr1} that for all
large enough $i$'s:
\begin{equation}\label{central1}
\Cste{31}\Prob{\Event{24}(i)}\leq \Prob{\Event{22}(i)} \leq
\Prob{\Event{24}(i)}.
\end{equation}
We now use Proposition \ref{estimonte} to check that:
\begin{equation}\label{central2}
\log\Pare{\Prob{\Event{24}(i)}} \underset{i\to\infty}{\sim}
-\frac{K^\#}{K}\frac{f(t_i)}{a^{-1}\Pare{e^i-e^{i/2}}}
\underset{i\to\infty}{\sim}-\frac{K^\#}{K}\log i,
\end{equation}
where we used the regular variation of $a(\cdot)$ for the last
equivalence. In particular, combining this with (\ref{central1}) and
the fact that $K^\#/K < 1$ show that
$\sum_{i}\P(\Event{22}(i))=\infty$. We now estimate
$\P(\Event{22}(i)\cap\Event{22}(j))$ for $i$ large enough and for $j>i$.
\begin{eqnarray*}
\Event{22}(i)\cap\Event{22}(j) &\subset&
\Event{24}(i)\cap\Event{24}(j)\\
&\subset&
\Event{24}(i)\cap\left\{\Pare{\theta_{f(t_i)/K}\V}^\#_{f(t_j)/K -
f(t_i)/K} \leq e^j - e^{j/2}\right\}.
\end{eqnarray*}
Hence, from the independence and the stationarity of the increments of $\V$
(at integer times), combined with Proposition \ref{Encadr2}, for all $j$
large enough (i.e. all $i$ large enough):
\begin{eqnarray*}
\Prob{\Event{22}(i)\cap\Event{22}(j)} &\leq&
\Prob{\Event{24}(i)}\Prob{\V^\#_{f(t_j)/K -
f(t_i)/K} \leq e^j - e^{j/2}}\\
&\leq&\Cste{32}\frac{\Prob{\Event{24}(i)}\Prob{\Event{24}(j)}}{\P\big(\V^\#_{f(t_i)/K}\leq
e^j - e^{j/2}\big)}.
\end{eqnarray*}
Using Lemma \ref{regvar}, one may check after a few lines of
calculus that for all $i$ sufficiently large, $\exp(j) -
\exp(j/2) \geq a^{-1}(f(t_i)/K)$ whenever $j-i \geq \log i$ thus
\begin{equation*}
\Prob{\V^\#_{f(t_i)/K}\leq e^j - e^{j/2}} \geq
\Prob{\frac{\V^\#_{f(t_i)/K}}{a\Pare{f(t_i)/K}}\leq 1}.
\end{equation*}
Since the r.h.s. of the last equation converges to $\P(\S^\#_1 \leq
1)\neq 0$ as $i$ goes to infinity we deduce that for all $i$ large
enough and all $j-i \geq \log i$:
\begin{equation*}
\Prob{\V^\#_{f(t_i)/K}\leq e^j - e^{j/2}} \geq \Cste{33} >0
\end{equation*}
Finally, for all $i$ large enough and for all $j\geq i$:
\begin{equation}\label{central3}
\Prob{\Event{22}(i)\cap\Event{22}(j)} \leq
\left\{%
\begin{array}{ll}
    \Prob{\Event{22}(i)} & \hbox{if $0 \leq j-i < \log i$}, \\
    \Cste{34}\Prob{\Event{22}(i)}\Prob{\Event{24}(j)} & \hbox{if $j-i \geq \log i$}. \\
\end{array}%
\right.
\end{equation}
Combining (\ref{central1}),(\ref{central2}) and (\ref{central3}), we see that
\begin{equation*}
\liminf_{n\rightarrow\infty}\sum_{i,j\leq n}
\Prob{\Event{22}(i)\cap\Event{22}(j)}\Big/\Big(\sum_{i\leq n}\Prob{\Event{22}(i)}\Big)^2
\leq \Cste{35},
\end{equation*}
thus the Borel-Cantelli Lemma of
\cite{Koche1} yields $\Prob{\Event{22}(i)\hbox{ i.o.}} >
1/\Cste{35}$. We now use a classical 0-1 argument (compare with
\cite{HuShi1}, p1511 for details) to conclude that
$\Prob{\Event{22}(i)\hbox{ i.o.}} = 1$. This proved that, with
probability $1$:
\begin{equation*}
\limsup_{t\rightarrow\infty}\frac{\overline{X}_t}{f(t)} \geq
\frac{1}{K^\#}.
\end{equation*}
Moreover, the value of $K^\#$ when the process $\V$ is completely
asymmetric case was calculated in Corollary \ref{valeurKdiese}.
\end{proof}

\subsection{Proof of Theorem \ref{MainTheo2}}

\begin{lemm}\label{lemm1theo2}
For all $t$ large enough and all $\frac{a^{-1}(\log t)}{(\log\log)^{2/q}}\leq v \leq a^{-1}(\log t)$:
\begin{equation*}
\Prob{\overline{X_t}< v} \leq
\Cste{36}\frac{b^{-1}(v)}{b^{-1}(a^{-1}(\log t))}.
\end{equation*}
\end{lemm}

\begin{proof} In the following, we assume that $t$
is a very large number (thus $v$ is also large). From
(\ref{egaliteset}) and Proposition \ref{propI1} and \ref{propI2}, we
deduce:
\begin{eqnarray*}
\Prob{\overline{X}_t < v} &\leq &\Prob{I_1(v) \geq \frac{t}{2}} + \Prob{I_2(v) \geq \frac{t}{2}} \\
&\leq & \Prob{\V^\#_v \geq \log\frac{t}{2} - (\log v)^4}\\
&&+ \Prob{\widetilde{U}_{\V}(\overline{\V}_v + (\log v)^4) \geq \log \frac{t}{2}} + \Cste{37}e^{-(\log v)^2}.
\end{eqnarray*}
Remind that $b(\cdot)$ is regularly varying with index $q<1$,
therefore using Corollary \ref{lemmefluc1} and Lemma \ref{regvar} we
find:
\begin{eqnarray*}
\Prob{\V^\#_v \geq \log\frac{t}{2} - (\log v)^4} &\leq& \Prob{\V^\#_v
\geq \frac{1}{2}\log t}\\
&\leq& \Cste{38}\frac{v}{a^{-1}\Pare{\log t}}\\
&\leq& \Cste{39}\frac{b^{-1}(v)}{b^{-1}(a^{-1}(\log t))}.
\end{eqnarray*}
It is also easy to check from the bounds on $v$ and the regular variation
of $a(\cdot)$ and $b(\cdot)$ that $$e^{-(\log v)^2}\leq
\frac{b^{-1}(v)}{b^{-1}\Pare{a^{-1}(\log t)}}.$$
We still have to prove a similar
bound for $\P(\widetilde{U}_{\V}(\overline{\V}_v + (\log v)^4)\geq
\log(t/2))$. Notice that for $b>a>0$, $\{\widetilde{U}_{\V}(a) \geq
b\} = \widetilde{\Lambda}'(a,b-a)$ hence using Proposition
\ref{tempssortie} and the independence of $(\V_x)_{x\geq 0}$ and
$(\V_{-x})_{x\geq 0}$:
\begin{eqnarray}
\nonumber&&\Prob{\widetilde{U}_{\V}(\overline{\V}_v + (\log v)^4) \geq
\log\frac{t}{2}} \\
\nonumber&&\qquad\leq \Cste{40}\Esp
{\frac{b^{-1}(a^{-1}(\overline{\V}_v + (\log v)^4))}{b^{-1}(a^{-1}(\log \frac{t}{2}))}} \\
\label{encoreunetd}&&\qquad =
\Cste{40}\frac{b^{-1}(v)}{b^{-1}(a^{-1}(\log\frac{t}{2}))}\Esp
{\frac{b^{-1}(a^{-1}(\overline{\V}_v + (\log
v)^4))}{b^{-1}(a^{-1}(a(v)))}}.
\end{eqnarray}
We now use Lemma \ref{regvar} for the regularly varying function 
$b^{-1}(a^{-1}(\cdot))$ to check that (\ref{encoreunetd}) is smaller than
\begin{equation*}
\Cste{41,\varepsilon}
\frac{b^{-1}(v)}{b^{-1}(a^{-1}(\log\frac{t}{2}))} \Esp
{\Pare{\frac{\overline{\V}_v + (\log v)^4}{a(v)}}^{\alpha q +
\varepsilon} + 1}.
\end{equation*}
Finally, since $q<1$, we can choose $\varepsilon$ small enough such that $\alpha
q + \varepsilon < \alpha$, therefore Corollary \ref{cormomentsup}
implies
\begin{equation*}
\Esp {\Pare{\frac{\overline{\V}_v + (\log v)^4}{a(v)}}^{\alpha q +
\varepsilon}} \leq \Esp
{\Pare{\frac{\overline{\V}_v}{a(v)}+1}^{\alpha q + \varepsilon}}
\leq \Cste{42,\varepsilon},
\end{equation*}
we conclude the proof noticing that $b^{-1}(a^{-1}(\log\frac{t}{2})) \sim b^{-1}(a^{-1}(\log t))$.
\end{proof}

\begin{lemm} \label{lemm2theo2}
For all $v$ large enough and for all $t>0$ we have
\begin{equation*}
\left\{ \overline{X}_t < v\right\} \supset
\left\{\widetilde{U}_{\V}(a(v)) \geq \log t\right\}\cap
\left\{\overline{\V}_{v/2} \geq 2a(v)\right\} \cap\Event{9}(v)
\end{equation*}
where $\Event{9}(v)$ was defined in Proposition \ref{propI2}  and satisfies 
$$\P(\Event{9}(v)^c)\leq \Cste{26}e^{-(\log v)^2}.$$
\end{lemm}
\begin{proof} Recall that relation (\ref{egaliteset}) gives 
$\{\overline{X}_t < v\} = \{I_1(v) + I_2(v) > t\}$ and notice
that $I_1(v) >0$ for all $v>0$ thus $\{\overline{X}_t < v\} \supset \{I_2(v) \geq t\}$.
We now use Proposition \ref{propI2} to see that for all $v$ large enough,
the event $\{\overline{X}_t < v\}$ contains
\begin{eqnarray*}
&&
\left\{\widetilde{U}_{\V}(\overline{\V}_{v-\frac{1}{2}}-(\log
v)^4)\geq \log t\right\}
\cap \left\{ \overline{\V}_{v-\frac{1}{2}} > (\log v)^{4} \right\}\cap\Event{9}(v)\\
&&\qquad\supset \left\{\widetilde{U}_{\V}(\overline{\V}_{v/2}-a(v))\geq
\log t\right\}
\cap \left\{ \overline{\V}_{v/2} \geq  2a(v) \right\}\cap\Event{9}(v)\\
&&\qquad\supset \left\{\widetilde{U}_{\V}(a(v))\geq \log t \right\} \cap
\left\{ \overline{\V}_{v/2} \geq  2a(v) \right\}\cap\Event{9}(v),
\end{eqnarray*}
where we used the fact that $x\mapsto\widetilde{U}_{\V}(x)$ is a non-decreasing
function and trivial inequalities $\overline{\V}_{v/2}
\leq {\V}_{v-1/2}$ and $(\log v)^4 \leq a(v)$ which hold for all
large enough $v$'s.
\end{proof}

\begin{proof}[Proof of Theorem \ref{MainTheo2}]
For any positive nondecreasing function $f$, let
\begin{equation*}
J(f) = \int^{\infty}\frac{b^{-1}(f(t))dt}{b^{-1}(a^{-1}(\log
t))t\log t}
\end{equation*}
(we do not specify the lower bound for the integral since we are only
concerned with the convergence of $J(f)$ at infinity). We easily
check using Lemma \ref{regvar} that $J(f) = \infty$ when
$f(t)=a^{-1}(\log t)/(\log \log t)^{1/(2q)}$ and that $J(f)<\infty$
when $f(t) = a^{-1}(\log t)/(\log \log t)^{2/q}$, therefore we may
assume without loss of generality that for $t$ large enough:
\begin{equation*}
\frac{a^{-1}(\log t)}{(\log\log t)^{2/q}}\leq f(t) \leq
\frac{a^{-1}(\log t)}{(\log\log t)^{1/(2q)}}.
\end{equation*}
We first assume that $J(f) < \infty$ and we define the sequence $t_i
= \exp(\exp i)$. Note that for $i$ large enough $a^{-1}(\log
t)/(\log\log t_i)^{2/q}\leq f(t_{i+1}) \leq a^{-1}(\log t_i)$ thus
we can use Proposition \ref{lemm1theo2}:
\begin{eqnarray*}
\Prob{\overline{X}_{t_i} < f(t_{i+1}) } &\leq&
\Cste{36}\frac{b^{-1}(f(t_{i+1}))}{b^{-1}(a^{-1}(\log t_i))}\\
&\leq&\Cste{43}\frac{b^{-1}(f(t_{i+1}))}{b^{-1}(a^{-1}(\log t_{i+2}))}\\
&\leq&\Cste{43}\int_{t_{i+1}}^{t_{i+2}}\frac{b^{-1}(f(t))dt}{b^{-1}(a^{-1}(\log
t)) t \log t},
\end{eqnarray*}
where we used that $b^{-1}(a^{-1}(\log t_{i+2}))\sim \exp(2\alpha
q)b^{-1}(a^{-1}(\log t_{i}))$ for the second inequality and the
monotonicity of $a^{-1}$,$b^{-1}$ and $f$ for the third inequality.
Since $J(f)<\infty$, we conclude that $\sum_{i}\P(\overline{X}_{t_i}
< f(t_{i+1}))<\infty$ and Borel-Cantelli's Lemma implies that
$\P(\overline{X}_{t_i} < f(t_{i+1})\hbox{ i.o.})=0$. For $t_i \leq
t\leq t_{i+1}$, we have $\overline{X}_t \geq \overline{X}_{t_i}$ and
$f(t_{i+1})\geq f(t_i)$ hence with probability $1$:
\begin{equation}\label{nuj}
\liminf_{t\rightarrow\infty}\frac{\overline{X}_t}{f(t)} \geq 1
\hbox{ a.s.}
\end{equation}
Changing $f$ for $C f$ for any $C>0$ does not alter the convergence
of $J(f)$ thus the $\liminf$ in (\ref{nuj}) is in fact infinite. We now
assume that $J(f)=\infty$. Using Lemma \ref{lemm2theo2}, for $i$
large enough:
\begin{equation*}
\left\{\overline{X}_{t_i} \leq f(t_i)\right\} \supset
\Event{9}(f(t_i))\cap \Event{26}(i),
\end{equation*}
where $\Event{26}(i)  =  \Event{27}(i)\cap\Event{28}(i)$ with
\begin{eqnarray*}
\Event{27}(i) & = & \{\widetilde{U}_{\V}(a(f(t_i))) \geq \log(t_i)\},\\
\Event{28}(i) & = & \{\overline{\V}_{f(t_i)/2} \geq 2a(f(t_i))\}.
\end{eqnarray*}
Since $\sum_{i}\P(\Event{9}(f(t_i))^c) < \infty$, it only remains
to prove that $\P(\Event{26}(i) \hbox{ i.o.})=1$. Results of
section \ref{subsectionD} imply that
$\lim_{i\rightarrow\infty}\P(\Event{28}(i)) = \P(\overline{\S}_{1/2}
\geq 2) >0$. Since $\Event{27}(i)$ and $\Event{28}(i)$ are independent
events, there exist a constant $\Cste{43}>0$ such that for all $i$
large enough:
\begin{equation}\label{llref1}
\Cste{43} \Prob{\Event{27}(i)}\leq \Prob{\Event{26}(i)} \leq
\Prob{\Event{27}(i)}.
\end{equation}
Notice that $a(f(t_i))$ and $\log(t_i)/a(f(t_i))$ both go to infinity
as $i$ goes to infinity. Using the estimate for the solution of the
exit problem (Proposition \ref{tempssortie}) and the regular variation of $b^{-1}\Pare{a^{-1}(\cdot)}$,
for all sufficiently large $i$'s:
\begin{equation}\label{llref2}
\Cste{44}\frac{b^{-1}(f(t_i))}{b^{-1}(a^{-1}(\log t_i))} \leq
\Prob{\Event{27}(i)}\leq
\Cste{45}\frac{b^{-1}(f(t_i))}{b^{-1}(a^{-1}(\log t_i))}.
\end{equation}
Combining the inequalities (\ref{llref1}) and (\ref{llref2}), the assumption that
$J(f)=\infty$  implies
$$\sum_{i}\P(\Event{26}(i)) = \infty.$$
We now
estimate $\P(\Event{26}(i)\cap\Event{26}(j))$. Let $g(i) =
\log(t_i)-a(f(t_i))$. It is easy to check that g is ultimately
increasing. Let us assume $i$ very large and let $j>i$. We can rewrite:
\begin{equation*}
\Event{27}(i)\cap\Event{27}(j) =\widetilde{\Lambda}'
\Pare{a(f(t_i)),g(i)}\cap\widetilde{\Lambda}'\Pare{a(f(t_j)),g(j)}.
\end{equation*}
There are two cases (which are not disjoint):
\begin{enumerate}
\item $(\V_{-n})_{n\geq 0}$ hits $(-\infty,-g(j)]$ before hitting
$[a(f(t_i)),\infty)$. We see from Proposition \ref{tempssortie} that
the probability of this case is less than
$\Cste{46}b^{-1}(f(t_i))/b^{-1}(a^{-1}(\log t_j))$.
\item $(\V_{-n})_{n\geq 0}$ hits $(-\infty,-g(i)]$ before hitting
$[a(f(t_i)),\infty)$ (i.e. $\Event{27}(i)$ happens) and the shifted
random walk $(\V_{-\widetilde{\sigma}_{\V}(a(f(t_i)))-n})_{n\geq 0}$ hits
$(-\infty,-g(j)]$ before hitting $[a(f(t_j)),+\infty)$ (the probability of
this event is smaller than $\P(\Event{27}(j))$). Using the Markov property
for the random walk $(\V_{-n})_{n\geq 0}$ we conclude that the probability 
of this case is smaller than $\P(\Event{27}(i))\P(\Event{27}(j))$.
\end{enumerate}
Combining (1) and (2) we deduce that $\P(\Event{27}(i)\cap\Event{27}(j))$ is smaller than
\begin{eqnarray*}
&&
\Prob{\Event{27}(i)}\Prob{\Event{27}(j)} +
\Cste{46}\frac{b^{-1}(f(t_i))}{b^{-1}(a^{-1}(\log t_j))}\\
&&\qquad \leq\Prob{\Event{27}(i)}\Prob{\Event{27}(j)} +
\frac{\Cste{46}}{\Cste{44}}\Prob{\Event{27}(i)}\frac{b^{-1}(a^{-1}(\log
t_i))}{b^{-1}(a^{-1}(\log t_j))},
\end{eqnarray*}
where we used (\ref{llref2}) for the second inequality.
Finally, using Lemma \ref{regvar} and (\ref{llref1}), we conclude
that for all $i$ large enough and all $j>i$:
\begin{eqnarray*}
\Prob{\Event{26}(i)\cap\Event{26}(j)} &\leq&\Prob{\Event{27}(i)\cap\Event{27}(j)}\\
&\leq& \Cste{47}\Pare{\Prob{\Event{26}(i)}\Prob{\Event{26}(j)} +
\Prob{\Event{26}(i)}e^{-\Cste{48}(j-i)}}
\end{eqnarray*}
hence
\begin{equation*}
\liminf_{n\rightarrow\infty}\sum_{i,j\leq n}
\Prob{\Event{26}(i)\cap\Event{26}(j)}\Big/\Big(\sum_{i\leq n}\Prob{\Event{26}(i)}\Big)^2
\leq \Cste{47}.
\end{equation*}
Just like for Theorem \ref{MainTheo1}, we apply the Borel-Cantelli
Lemma of \cite{Koche1} and a standard 0-1 argument to conclude that
$\P(\Event{26}(i) \hbox{ i.o.})=1$. Since this result still holds when changing
$f$ for $Cf$ for any $C>0$, we have proved that, with probability $1$,
\begin{equation*}
\liminf_{t\rightarrow\infty}\frac{\overline{X}_t}{f(t)}=0.
\end{equation*}
\end{proof}

\subsection{Proof of Theorem \ref{MainTheo2b}}
Just like the previous two theorems, the proof is based on the
following two lemmas.

\begin{lemm}
For all $t$ large enough and all $\lambda$ such that:
$$(\log\log t)^{1/4}\leq\lambda\leq (\log\log t)^4,$$
we have
$$\Prob{X^*_t < \frac{a^{-1}(\log t)}{\lambda}}\leq \frac{\Cste{49}}{\lambda^2}.$$
\end{lemm}

\begin{proof}
We use the notation $v = a^{-1}(\log t)/\lambda$.
According to (\ref{egaliteset}) we have:  $$\left\{\overline{X}_t <
v\right\} = \left\{I_1(v) + I_2(v) > t\right\}$$ where $I_1$ and
$I_2$ were defined in (\ref{expI1}) and (\ref{expI2}). Using a
symmetry argument:
$$\left\{\underline{X}_t > -v\right\} = \left\{\widetilde{I}_1(v) + \widetilde{I}_2(v) > t\right\},$$
where $\widetilde{I}_1$ and $\widetilde{I}_2$ are given again by the
formulas (\ref{expI1}) and (\ref{expI2}) by simply changing the
process $(\V_x)_{x\in\R}$ for $(\V_{-x})_{x\in\R}$. Combining these
equalities, we get:
\begin{equation}\label{pref4}
\Big\{  X^*_t < v  \Big\}  = \Big\{I_1(v) + I_2(v) > t\Big\} \cap
\left\{\widetilde{I}_1(v) + \widetilde{I}_2(v) > t\right\},
\end{equation}
hence $\P(X^*_t < v)$ is smaller than
\begin{equation*}
\Prob{I_1(v) + I_2(v) > t\hbox{ ,
}\overline{\V}_v \leq \overline{\V}_{-v}} + \Prob{
\widetilde{I}_1(v) + \widetilde{I}_2(v) > t \hbox{ ,
}\overline{\V}_{-v}\leq \overline{\V}_v}.
\end{equation*}
It is clear from a symetry argument that we only need to prove the bound for
the first member of the last equation. Notice that:
\begin{eqnarray}
&&\nonumber\Prob{I_1(v) + I_2(v) > t\hbox{ , }\overline{\V}_v \leq
\overline{\V}_{-v}}\\
\label{pref1}&&\qquad\leq\Prob{\frac{1}{4}\log t \leq \overline{\V}_v \leq \overline{\V}_{-v}}\\
\label{pref2}&&\qquad\qquad + \Prob{I_1(v) \geq \frac{t}{2} , \overline{\V}_v \leq \frac{\log t}{4}} \\
\label{pref3}&&\qquad\qquad + \Prob{I_2(v) \geq \frac{t}{2}
,\overline{\V}_{v}\leq \overline{\V}_{-v},\overline{\V}_v \leq
\frac{\log t}{4}}.
\end{eqnarray}
We deal with each term separately. First, using independence of
$(\V_x)_{x\geq 0}$ and $(\V_{-x})_{x\geq 0}$ we see that
(\ref{pref1}) is smaller than
$$
\Prob{\overline{\V}_v \geq \frac{1}{4}\log
t}\Prob{\overline{\V}_{-v} \geq \frac{1}{4}\log t} \leq
\frac{\Cste{49}}{\lambda^2},
$$
where we used Corollary \ref{lemmefluc1} for the last inequality. We
now turn our attention to (\ref{pref2}). Using Proposition
\ref{propI1}, we check that this probability is smaller than
$$
\Prob{\V^\#_v \geq \log \frac{t}{2} - \log^4 v\hbox{ ,
}\overline{\V}_v \leq \frac{1}{4}\log t} + \Cste{25}e^{-\log^2 v}.
$$
For $t$ large enough, using the Markov property:
\begin{eqnarray*}
&&\Prob{\V^\#_v \geq \log \frac{t}{2} - \log^4 v\hbox{ , }\overline{\V}_v
\leq \frac{\log t}{4}} \\
&&\qquad\leq\Prob{\V^\#_v \geq \frac{\log t}{2}\hbox{ , }
\overline{\V}_v \leq \frac{\log t}{4}} \\
&&\qquad\leq \Prob{\sigma_{\V}(-\frac{\log t}{4})\leq v ,
\Pare{\theta_{\sigma_{\V}(-\frac{\log t}{4})}\V}^\#_v \geq \frac{\log t}{2}}\\
&&\qquad\leq\Prob{\underline{\V}_{v} \leq -\frac{\log t}{4}}\Prob{\V^\#_v\geq \frac{\log t}{2}}\\
&&\qquad\leq\frac{\Cste{50}}{\lambda^2},
\end{eqnarray*}
where we used again Corollary \ref{lemmefluc1} for the last line.
Note also that from the bound on $\lambda$, we have $e^{-\log^2 v}
\leq 1/\lambda^2$ for all $t$ large enough. This gives the desired
bound for (\ref{pref2}). It remains to prove the existence of a similar bound for
(\ref{pref3}). We first use Proposition \ref{propI2} to see that, for
all $t$ large enough, (\ref{pref3}) is smaller than
$$
\P\Big(\widetilde{U}_{\V}(\overline{\V}_v + \log^4 v) \geq
\log\frac{t}{2}\hbox{ , } \overline{\V}_{v}\leq \overline{\V}_{-v}
\hbox{ , } \overline{\V}_v \leq \frac{1}{4}\log t\Big) +
\Cste{26}e^{-\log^2 v}.
$$
We can rewrite:
\begin{eqnarray*}
&&\Big\{\widetilde{U}_{\V}(\overline{\V}_v + \log^4 v) \geq
\log\frac{t}{2}\hbox{ , }
\overline{\V}_{v}\leq \overline{\V}_{-v} \hbox{ , } \overline{\V}_v \leq \frac{1}{4}\log t\Big\}\\
&&\qquad= \Big\{\widetilde{\sigma}_{\V}\Big(\overline{\V}_v + \log^4 v - \log\frac{t}{2}\Big) <
\widetilde{\sigma}_{\V}\Pare{\overline{\V}_v+\log^4 v} \hbox{ , }\\
&&\qquad\qquad\widetilde{\sigma}_{\V}(\overline{\V}_v) \leq v\hbox{ , }
\overline{\V}_v \leq \frac{1}{4}\log t\Big\}\\
&&\qquad\subset \Big\{\widetilde{\sigma}_{\V}\Big(- \frac{\log t}{2}\Big) <
\widetilde{\sigma}_{\V}\Pare{\overline{\V}_v+\log^4 v} \hbox{ , }
\widetilde{\sigma}_{\V}(\overline{\V}_v) \leq v\Big\}\\
&&\qquad\subset \Big\{\widetilde{\sigma}_{\V}\Big(- \frac{\log t}{2}\Big) <
\widetilde{\sigma}_{\V}(\overline{\V}_v) \leq v\Big\}\\
&&\qquad\qquad\cup
\Big\{\widetilde{\sigma}_{\V}(\overline{\V}_v) <
\widetilde{\sigma}_{\V} \Big(- \frac{\log t}{2}\Big) <
\widetilde{\sigma}_{\V}\Pare{\overline{\V}_v+\log^4 v}\Big\}.
\end{eqnarray*}
Notice that on the event $\{\widetilde{\sigma}_{\V}(-(\log t)/2) <
\widetilde{\sigma}_{\V}(\overline{\V}_v) \leq v\}$, the process
$(\V_{-x})_{x\geq 0}$ hits $(-\infty,-(\log t/2)]$ before time $v$
and from this time on it hits $[0,\infty)$, again before time $v$, hence the
Markov property with the stopping time
$\widetilde{\sigma}_{\V}(-(\log t)/ 2)$ and Corollary
\ref{lemmefluc1} yields:
\begin{eqnarray*}
&&\Prob{\widetilde{\sigma}_{\V}\Pare{- \frac{\log t}{2}} <
\widetilde{\sigma}_{\V}(\overline{\V}_v) \leq v}\\
&&\qquad\leq
\Prob{\underline{\V}_{-v} \leq -\frac{\log
t}{2}}\Prob{\overline{\V}_{-v} \geq \frac{\log t}{2}}\leq
\frac{\Cste{51}}{\lambda^2}.
\end{eqnarray*}
It is also easy to check from the Markov property of
$(\V_{-x})_{x\geq 0}$ applied with the stopping time
$\widetilde{\sigma}_{\V}(\overline{\V}_v)$ that the probability of
the event $\{\widetilde{\sigma}_{\V}(\overline{\V}_v) <
\widetilde{\sigma}_{\V}(-(\log t)/2) <
\widetilde{\sigma}_{\V}\Pare{\overline{\V}_v+\log^4 v}\}$ is smaller
than the probability that the random walk $(\V_{-x})_{x\geq 0}$ hits
$(-\infty,-(\log t)/2]$ before it hits $[\log^4 v,\infty)$. Using the
estimate for the exit problem (Proposition \ref{tempssortie}) and the
regular variation of $b\Pare{a^{-1}(\cdot)}$, for $t$ large enough, we have:
\begin{eqnarray*}
&&\Prob{\widetilde{\sigma}_{\V}(\overline{\V}_v) <
\widetilde{\sigma}_{\V} \Pare{- \frac{\log t}{2}} <
\widetilde{\sigma}_{\V} \Pare{\overline{\V}_v+\log^4 v}}\\
&&\qquad\leq\Cste{52} \frac{b^{-1}\Pare{a^{-1}\Pare{(\log
v)^4}}}{b^{-1}\Pare{a^{-1}\Pare{\frac{\log t}{2}}}}\leq
\frac{1}{\lambda^2},
\end{eqnarray*}
so we conclude that (\ref{pref3}) is smaller than
$\Cste{53}/\lambda^2$.
\end{proof}

\begin{lemm}
for all $t$ large enough and all $(\log\log t)^{1/4}\leq\lambda\leq (\log\log t)^4$
we have:
$$
\left\{X^*_t \leq \frac{a^{-1}(\log t)}{\lambda}\right\}\supset
\left\{ \V_{v-\frac{1}{2}} \geq 2\log t \hbox{ , }
\V_{-v+\frac{1}{2}} \geq 2\log t\right\}\cap\Event{29}(v),
$$
where $v = a^{-1}(\log t)/\lambda$  and where
$$
\Prob{\Event{29}^c} \leq e^{-\Cste{54}\lambda^{1/4}}.
$$
\end{lemm}

\begin{proof}
Recall the definition for $\widetilde{I}_1$ given in the last lemma.
We assume $t$ very large. From (\ref{pref4}), we get
\begin{equation}
\Big\{  X^*_t < v  \Big\}  \supset \Big\{I_1(v) > t\Big\} \cap
\left\{\widetilde{I}_1(v) > t\right\},
\end{equation}
and Proposition \ref{propI1} yields
$$\Big\{I_1(v) > t\Big\} \supset
\Big\{\V^\#_{v-1/2}
> \log t + \log^4 v\Big\}\cap \Event{8}(v),$$
with $\P(\Event{8}^c)\leq \Cste{25}\exp(-\log^2 v)$. Similarly,
since $\widetilde{I}_1$ is obtained just like $I_1$ by changing
$(\V_x)_{x\in\R}$ for $(\V_{-x})_{x\in\R}$ in (\ref{expI1}), we also have
$$\{\widetilde{I}_1(v) > t\} \supset \{\V^\#_{-v+1/2} > \log t +
\log^4 v\}\cap \Event{30}(v)$$
where $\Event{30}(v)$ is a measurable set such that
$\P(\Event{30}^c)\leq \Cste{55}\exp(-\log^2 v)$. Let us define the event $\Event{29} =
\Event{8}\cap\Event{30}$. One may check from the bounds on $\lambda$
that $\Prob{\Event{29}^c} \leq \exp\Pare{-\Cste{54}\lambda^{1/4}}$ and
\begin{eqnarray*}
&&\big\{  X^*_t < v  \big\}\\
&&\qquad\supset \big\{\V^\#_{v-1/2} > \log t + \log^4 v\big\}\cap
\big\{\V^\#_{-v+1/2} > \log t + \log^4 v\big\}\cap \Event{29}(v)\\
&&\qquad\supset \big\{\V^\#_{v-1/2} \geq 2\log t\big\}\cap \big\{\V^\#_{-v+1/2} \geq 2\log t\big\}\cap \Event{29}(v)\\
&&\qquad\supset \big\{\V_{v-1/2} \geq 2\log t\big\}\cap \big\{\V_{-v+1/2} \geq 2\log
t\big\}\cap \Event{29}(v).
\end{eqnarray*}
\end{proof}

\begin{proof}[Proof of Theorem \ref{MainTheo2b}.]
This theorem is an easy consequence (using similar technics as in the
proof of Theorem \ref{MainTheo2}) of the last two lemmas and of
Proposition \ref{precislemmefluc}. We feel free to omit it.
\end{proof}

\subsection{Proof of Theorem \ref{MainTheo3}}
\begin{prop} \label{lastone} We have:
\begin{equation*}
\frac{1}{a(v)}\Pare{\log \sigma_{X}(v) - \V^\#_v \vee
\widetilde{U}_{\V}(\overline{\V}_v)} \overset{\hbox{
Prob.}}{\underset{v\to\infty}{\longrightarrow}} 0.
\end{equation*}
\end{prop}
The proof of this Proposition  is very similar to that of
Proposition $11.1$ of \cite{HuShi1} using the estimates for $I_1$
and $I_2$ obtained in Propositions \ref{propI1} and \ref{propI2}, we
therefore skip the details.

\begin{proof}[Proof of Theorem \ref{MainTheo3}]
Let $\lambda>0$ and let $v$ be a large number:
\begin{eqnarray*}
\P\Pare{\frac{\overline{X}_v}{a^{-1}(\log v)} \geq \lambda}
&=& \P\Pare{\log \sigma_{X}(\lambda a^{-1}(\log v))\leq \log v}\\
&=& \P\Pare{\frac{\log \sigma_{X}(x)}{c(x)} \leq
\frac{1}{\lambda^{1/\alpha}}},
\end{eqnarray*}
with the change of variable $x=\lambda a^{-1}(\log v)$ and where
\begin{equation}\label{laxxlls}
c(x) = \frac{\lambda^{1/\alpha}a(x/\lambda)}{a(x)} \underset{x\to\infty}{\sim} a(x)
\end{equation}
Results of section \ref{subsectionD} insure that
$(\V^\#_x \vee\widetilde{U}_{\V}(\overline{\V}_x))/a(x)$ converges in
law as $x\to\infty$ towards $\S^\#_1
\vee\widetilde{U}_{\S}(\overline{\S}_1)$ whose cumulative function
is continuous, hence it follows from Proposition \ref{lastone} and
from (\ref{laxxlls}) that
\begin{equation*}
\lim_{v\rightarrow\infty} \P\Pare{\frac{\overline{X}_v}{a^{-1}(\log
v)} \geq \lambda} = \Prob{ \S^\#_1
\vee\widetilde{U}_{\S}(\overline{\S}_1) \leq
\frac{1}{\lambda^{1/\alpha}}}.
\end{equation*}
This proves the convergence in law of $\overline{X}_v/a^{-1}(\log
v)$ towards the non degenerate random variable $\Xi = (\S^\#_1
\vee\widetilde{U}_{\S}(\overline{\S}_1))^{-\alpha}$ as $v\to\infty$.
Let us calculate the Laplace transform of this law when $\S$ is
completely asymmetric. Recall the notation $\tau^\#_x$ and $\tau_x$
defined in section \ref{subsectionA}. Let also $r_1$ be the stopping
time:
\begin{equation*}
r_1 = \inf\Pare{x\geq 0\hbox{ , }(\S_{-t})_{t\geq 0} \hbox{ hits }
(-\infty , -(1-x)) \hbox{ before it hits } (x,\infty) }.
\end{equation*}
Using the scaling property of $\S$:
\begin{eqnarray*}
\Prob{(\S^\#_1 \vee\widetilde{U}_{\S}(\overline{\S}_1))^{-\alpha}
\leq \lambda} &=&
\Prob{\S^\#_\lambda \vee\widetilde{U}_{\S}(\overline{\S}_\lambda) \geq 1}\\
&=&\Prob{ \tau^\#_{1}\wedge\tau_{r_1}\leq \lambda},
\end{eqnarray*}
therefore $\Xi$ and $\tau^\#_{1}\wedge\tau_{r_1}$ have the same law. Let
us first assume that $\S$ has no positive jumps and recall that
$(-\S_{-t}\hbox{ , }t\geq 0)$ and $(\S_{t}\hbox{ , }t\geq 0)$ have the
same law. It follows from the well known solution of the exit problem
for a completely asymmetric Levy process via its scale function $W$
(c.f. \cite{Berto1} , p194) that:
\begin{eqnarray*}
\Prob{r_1 > x} &=& \Prob{\hbox{$(\S_{-t})_{t\geq 0}$ hits $(x,\infty)$ before it hits $(-\infty,-(1-x))$}}\\
&=& 1 - \Prob{\hbox{$(\S_{t})_{t\geq 0}$ hits $(1-x,\infty)$ before it hits $(-\infty,-x)$}}\\
&=& 1 - \frac{W(x)}{W(1)},
\end{eqnarray*}
and it is known that in our case $W(x) = x^{\alpha-1}/
\Gamma(\alpha)$, hence the density of $r_1$ is
\begin{equation*}
\Prob{r_1 = dx} = \frac{\alpha-1}{x^{2-\alpha}}dx \hbox{ for
$x\in(0,1)$}.
\end{equation*}
Using Proposition \ref{pourloi} and the independence of
$(\S_{t})_{t\geq 0}$ and $(\S_{-t})_{t\geq 0}$  we have for $q\geq
0$:
\begin{eqnarray*}
\Esp{e^{-q \tau^\#_1\wedge\tau_{r_1}}} &=&
\int_{0}^{1}\Esp{e^{-q\tau^\#_{1}\wedge\tau_{x}}}\frac{\alpha-1}{x^{2-\alpha}}dx\\
& = & \frac{\alpha -1}{\Mittag_\alpha(q)}\int_{0}^{1}
\frac{\Mittag_{\alpha}(q(1-x)^\alpha)}{x^{2-\alpha}}dx\\
& = & \frac{\alpha
-1}{\Mittag_\alpha(q)}\sum_{n=0}^{\infty}\frac{q^n}{\Gamma(1+\alpha
n)}\int_{0}^{1}\frac{(1-x)^{\alpha n}}{x^{2-\alpha}}dx,
\end{eqnarray*}
but
\begin{equation*}
\frac{1}{\Gamma(1+\alpha n)}\int_{0}^{1}\frac{(1-x)^{\alpha
n}}{x^{2-\alpha}}dx = \frac{\Gamma(\alpha-1)}{\Gamma(\alpha(n+1))},
\end{equation*}
hence
\begin{eqnarray*}
\Esp{e^{-q \tau^\#_1\wedge\tau_{r_1}}} &=&  \frac{\Gamma(\alpha)}{\Mittag_\alpha(q)}
\sum_{n=0}^{\infty}\frac{q^n}{\Gamma(\alpha(n+1))}\\
& = & \Gamma(\alpha+1)
\frac{\Mittag_\alpha'(q)}{\Mittag_\alpha'(q)}.
\end{eqnarray*}
We now assume that $\S$ has no negative jumps. Just like in the
previous case, we can calculate the density of $r_1$ from the scale
function and we find $\P(r_1 = dx) = (\alpha -1)/(1-x)^{2-\alpha}$
for $x\in(0,1)$ thus using Proposition \ref{pourloi}:
\begin{eqnarray*}
\Esp{e^{-q \tau^\#_1\wedge\tau_{r_1}}} &=& \int_{0}^{1}
\Esp{e^{-q\tau^\#_{1}\wedge\tau_{x}}}\frac{\alpha-1}{x^{2-\alpha}}dx\\
&=& (\alpha-1)\int_{0}^{1}\frac{\Mittag_\alpha(qx^\alpha)}{(1-x)^{2-\alpha}}dx\\
&& - \frac{\Mittag_\alpha'(q)\alpha(\alpha-1) q}{\alpha q
\Mittag_\alpha''(q) + (\alpha
-1)\Mittag_\alpha'(q)}\int_{0}^{1}\frac{x^{\alpha-1}\Mittag_\alpha'(qx^\alpha)}{(1-x)^{2-\alpha}}dx.
\end{eqnarray*}
We already calculated the first integral:
\begin{equation*}
\int_{0}^{1}\frac{\Mittag_\alpha(qx^\alpha)}{(1-x)^{2-\alpha}}dx =
\int_{0}^{1}\frac{\Mittag_\alpha(q(1-y)^\alpha)}{y^{2-\alpha}}dy =
\frac{\Gamma(\alpha+1)}{\alpha-1}\Mittag_\alpha'(q).
\end{equation*}
As for the second integral:
\begin{equation*}
\int_{0}^{1}\frac{x^{\alpha-1}\Mittag_\alpha'(qx^\alpha)}{(1-x)^{2-\alpha}}dx =
\sum_{n=0}^{\infty}\frac{(n+1)q^n}{\Gamma(\alpha(n+1)+1)}
\int_{0}^{1}\frac{x^{\alpha(n+1)-1}}{(1-x)^{2-\alpha}}dx,
\end{equation*}
and it is known that
\begin{equation*}
\int_{0}^{1}\frac{x^{\alpha(n+1)-1}}{(1-x)^{2-\alpha}}dx =
\frac{\Gamma(\alpha(n+1))\Gamma(\alpha-1)}{\Gamma(\alpha(n+2)-1)},
\end{equation*}
hence
\begin{eqnarray*}
&&\int_{0}^{1}\frac{x^{\alpha-1}\Mittag_\alpha'(qx^\alpha)}{(1-x)^{2-\alpha}}dx\\
&&\quad=
\frac{\Gamma(\alpha-1)}{\alpha}\sum_{n=0}^{\infty}\frac{q^n}{\Gamma(\alpha(n+2)-1)}\\
&&\quad= \Gamma(\alpha-1)\sum_{n=0}^{\infty}\frac{(n+2)(\alpha(n+2)-1)q^n}{\Gamma(\alpha(n+2)+1)}\\
&&\quad= \Gamma(\alpha-1)\Pare{\alpha\sum_{n=0}^{\infty}\frac{(n+1)(n+2)q^n}{\Gamma(\alpha(n+2)+1)}
+ (\alpha-1)\sum_{n=0}^{\infty}\frac{(n+2)q^n}{\Gamma(\alpha(n+2)+1)}}\\
&&\quad= \frac{\Gamma(\alpha-1)}{q}\Pare{q\alpha\Mittag_\alpha''(q) +
(\alpha-1)\Mittag_\alpha'(q) - \frac{\alpha-1}{\Gamma(\alpha+1)}}.
\end{eqnarray*}
Putting the pieces together, we conclude:
\begin{equation*}
\Esp{e^{-q \tau^\#_1\wedge\tau_{r_1}}} =
\frac{(\alpha-1)\Mittag_\alpha'(q)}{\alpha q \Mittag_\alpha''(q) +
(\alpha-1)\Mittag_\alpha'(q)}.
\end{equation*}
\end{proof}

\section{Comments}
\subsection{The case where $\V$ is a stable process.}
In the whole paper, we assumed $\V$ to be a random walk in the
domain of attraction of a stable process $\S$. Let us now
assume that $\V$ itself is a strictly stable process (such that
$|\V|$ is not a subordinator) and let us explain why Theorems
$\ref{MainTheo1}-\ref{MainTheo3}$ still hold in this case.
It is clear that all the results dealing with the fluctuations
of $\V$ remain unchanged (in fact, they even take a nicer form
since we can now choose $a(x)=x^{\alpha}$ and $b(x)=x^q$). Notice
also that we did not use the fact that $\V$ was a random walk in
the proofs of the theorems in section $4$. Indeed, the only time 
we really used the assumption that $\V$ was flat on the intervals
$(n,n+1)\hbox{ , }n\in\Z$ was in the proofs of Propositions \ref{propI1}
and \ref{propI2} (we needed to make sure that $\V$ spends ``enough'' time
around its local extremas). Looking closely at those two proofs, we see 
that they will still hold if we can show that there exist a measurable
event $\Event{31}(v)$ such that:
\begin{enumerate}
\item there exists $\Cste{56}$ such that $\Prob{\Event{31}(v)^c}\leq \Cste{56}\exp(-\log^2 v)$.
\item On $\Event{31}(v)$, any path of $\V$ is such that for all 
$x\in[-\widetilde{\sigma}_{\V}(\overline{\V}_v+\log^4 v) , v]$, we have
$|\V_{y} -\V_{x}|\leq 1$ either for all $y$ in $[x,x+\exp(-\log^3 v)]$ or
for all $y$ in the interval $[x-\exp(-\log^3 v),x]$.
\end{enumerate}
Let us quickly explain how we can construct this event. Define the sequence of
random variables $(\gamma_n)_{n\in\Z}$:
\begin{equation*}
\left\{%
\begin{array}{l}
\gamma_0 = 0, \\
\gamma_{n+1}  = \inf( t>\gamma_n \hbox{ , } | \V_{t}-\V_{\gamma_n}|\geq\frac{1}{2}) \hbox{ for $n\geq0$,}\\
\gamma_{-n-1}  = \inf(t<-\gamma_n \hbox{ , } |\V_{t}-\V_{-\gamma_n}|\geq\frac{1}{2}) \hbox{ for $n\geq0$.}
\end{array}%
\right.
\end{equation*}
Let us set
\begin{eqnarray*}
\Event{32}(v) &=& \left\{ \gamma_{i+1} - \gamma_i > 2e^{-\log^3 v}
\hbox{ for all $-e^{\frac{1}{2}\log^{3} v} \leq i\leq e^{\frac{1}{2}\log^{3}v}$}\right\},\\
\Event{33}(v) &=& \left\{ \gamma_{-[e^{\frac{1}{2}\log^{3} v}]}> e^{\log^{5/2}v}
\hbox{ , }  \gamma_{[e^{\frac{1}{2}\log^{3} v}]}> e^{\log^{5/2}v}\right\},\\
\Event{34}(v) &=& \left\{ \widetilde{\sigma}_{\V}(\overline{\V}_v+\log^4 v) \leq e^{\log^{5/2} v}\right\},\\
\Event{31}(v) &=& \Event{32}(v)\cap\Event{33}(v)\cap\Event{34}(v).
\end{eqnarray*}
It is clear that condition ($2$) holds for $\Event{31}$. We now assume that $v$ is very large. We have:
\begin{equation*}
\Prob{\Event{32}(v)^c}\leq 2e^{\frac{1}{2}\log^3 v}\Prob{\gamma_1 \leq 2e^{-\log^3 v}}\leq \Cste{57}e^{-\frac{1}{2}\log^3 v},
\end{equation*}
where we used the relation $\P(\gamma_1\leq x) = \P(\V^*_{x}\geq \frac{1}{2})$
and Corollary \ref{lemmefluc1} for the last inequality. Using Cramer's large deviation theorem,
it is easy to check that $\Prob{\Event{33}(v)^c} \leq e^{-v}$ (in fact, we can obtain a much
better bound). We also have $\P(\Event{34}(v)^c) \leq 3 e^{-\log^2 v}$ (compare with  the
proof page \pageref{al89} of the inequality $\P(\Event{16}(v)^c) \leq 3e^{-\log^2 v}$ for
details). Thus condition ($1$) holds.

\subsection{Non-symetric environments.} In the whole paper, in order to avoid even more complicated
notations, we assumed that the processes $(\V_x \hbox{ , }x\geq 0)$ and $(-\V_{-x}\hbox{ ,
}x\leq 0)$ have the same law. However it is easy to see that this assumption can be
relaxed. Indeed, we may swap assumption \ref{hyp1} for the following:
\begin{hypo} \label{hyp2} $(\V_n)_{n\geq 0}$ and $(\V_{-n})_{n\geq 0}$
are independent random walks and there exists a positive  sequence
$\Pare{a_n}_{n\geq 0}$ such that
\begin{equation*}
\frac{\V_n}{a_n}\underset{n\to\infty}{\overset{\hbox{law}}{\longrightarrow}}\S^1\hbox{
and }
\frac{-\V_{-n}}{a_n}\underset{n\to\infty}{\overset{\hbox{law}}{\longrightarrow}}\S^2,
\end{equation*}
where $\S^1$ and $\S^2$ are random variables whose law are strictly
stable with respective parameters $(\alpha,p_1)$ and $(\alpha,p_2)$
and whose densities are everywhere positive on $\R$.
\end{hypo}
It is crucial to assume that the norming sequence $(a_n)$ may be
chosen to be the same for both random walk (in order to keep the
results of functional convergence of section \ref{subsectionD}) but
the positivity parameters $p_1$ and $p_2$ need not be the same.
Theorem \ref{MainTheo1}-\ref{MainTheo3} must be adapted in
consequences. For example, Theorem \ref{MainTheo1} now take the form:
\begin{theop} Under the annealed probability
$\P$, almost surely:
\begin{equation*}
\limsup_{t\to\infty}\frac{X_t}{a^{-1}\Pare{\log t}\log\log\log t} =
\frac{1}{K^{\#,1}},
\end{equation*}
where $K^{\#,1}$ depends only on $\S^1$ and is given by the formula:
\begin{equation*}
K^{\#,1} = -\lim_{t\to\infty}\frac{1}{t}\log \Prob{\sup_{0\leq u\leq
v\leq t}\Pare{\S^1_v - \S^1_u} \leq 1}.
\end{equation*}
Furthermore, when $\S^1$ is completely asymmetric:
$K^{\#,1}$ is given by:
\begin{equation*}
K^{\#,1} =
\left\{%
\begin{array}{l}
    \rho_1 (\alpha) \hbox{ when $\S^1$ has no positive jumps}, \\
    \rho_2 (\alpha) \hbox{ when $\S^1$ has no negative jumps}. \\
\end{array}%
\right.
\end{equation*}
\end{theop}
Let now $(\T_n)$ stands for the sequence of strict ascending ladder
index of the random walk $(\V_{-x})_{x\geq 0}$:
\begin{equation*}
\left\{%
\begin{array}{l}
    \T_0 =0, \\
    \T_{n+1} = \min\Pare{k>\T_n \hbox{ , } \V_{-k} > \V_{-\T_n} }.
\end{array}%
\right.
\end{equation*}
hence $\T_1$ is in the domain of attraction of a positive stable law
with index $p_2$ and we choose $b(\cdot)$ to be a continuous
positive increasing function such that $(b(n))_{n\geq 1}$ is a norming
sequence for $\T_1$. Theorem \ref{MainTheo2} now takes the form:
\begin{theop}  For any non decreasing function $f$ we have:
\begin{equation*}
\P\Pare{\sup_{s\leq t} X_s \leq f(t) \hbox{ i.o. }} =
\left\{%
\begin{array}{l}
    0 \\
    1
\end{array}%
\right. \Longleftrightarrow
\int^{\infty}\frac{b^{-1}\Pare{f(t)}dt}{b^{-1}\Pare{a^{-1}\Pare{\log
t}}t\log t}
\left\{%
\begin{array}{l}
    <\infty \\
    =\infty.
\end{array}%
\right.
\end{equation*}
In particular, with probability $1$:
\begin{equation*}
\liminf_{t\to\infty}\frac{\Pare{\log\log
t}^{\beta}}{a^{-1}\Pare{\log t}}\sup_{s\leq t}X_s =
\left\{%
\begin{array}{ll}
    0, & \hbox{if $\beta<1/p_2$}, \\
    \infty, & \hbox{if $\beta>1/p_2$}. \\
\end{array}%
\right.
\end{equation*}
\end{theop}
Theorems \ref{MainTheo2b} and \ref{MainTheo3} must be adapted
similarly. Note that for Theorem \ref{MainTheo3} we can again calculate
the Laplace transform of the limiting law when $\S^1$ and $\S^2$
have both completely asymmetric laws.

\subsection{Random walk in random environment.}
let us recall the connection between the diffusion in
random potential and the model of Sinai's random walk in random
environment. Let $\omega =(\omega_i)_{i\in\Z}$ be an i.i.d. family
of random variables in $(0,1)$ and define for each realization of
this family a Markov chain $(Z_n)_{n\geq 0}$ by $Z_0 = 0$ and
\begin{equation*}
\P\Pare{Z_{n+1} = Z_{n} + e\hbox{ | }Z_n = x , (\omega_i)_{i\in\Z}} =
\left\{%
\begin{array}{ll}
    \omega_x & \hbox{if $e=1$,}\\
    1-\omega_x & \hbox{if $e=-1$.}
\end{array}%
\right.
\end{equation*}
$(Z_n)$ is a random walk in the random environment $\omega$. We now
define the associated two-sided random walk $(\V_{n})_{n\in\Z}$ by
$\V_0=0$ and $\V_{n+1} - \V_n= \log\Pare{(1-\omega_n)/\omega_n}$ for
all $n\in\Z$. Let $X$ still denotes the random diffusion in the
random potential $\V$. The following result from Schumacher
\cite{Schum1} relates the two processes $X$ and $Z$:
\begin{prop}
Define the sequence $(\mu_{n})_{n\geq 0}$ by
\begin{equation*}
\left\{%
\begin{array}{l}
    \mu_0 =0,  \\
    \mu_{n+1} = \inf\Pare{t>\mu_n \hbox{ , } |X_{\mu_{n+1}} -
    X_{\mu_n}|=1}.
\end{array}%
\right.
\end{equation*}
Under the annealed probability $\P$, the sequence $(\mu_{n+1} -
\mu_n)_{n\geq 0}$ is i.i.d. and $\mu_1$ is distributed as the first
hitting time of $1$ of a reflected standard Brownian
motion. Moreover, for each realization of the environment $\omega$.
The processes $(X_{\mu_n})_{n\geq 0}$ and $(Z_n)_{n\geq 0}$ have
same law.
\end{prop}
Using this proposition, we can easily adapt Theorem
\ref{MainTheo1}-\ref{MainTheo3} for the random walk in random
environment $Z$ in the case where $\V_1 =
\log\Pare{(1-\omega_0)/\omega_0}$ verifies assumption \ref{hyp1} (see
section $10$ of \cite{HuShi1} for details). For example, Theorem
\ref{MainTheo2b} for $Z$ takes the form:
\begin{theop}
When $\S$ has jumps of both signs, we have
for any increasing positive sequence $(c_n)_{n\geq 0}$:
\begin{equation*}
\P\Pare{\sup_{k\leq n} |Z_k| \leq \frac{a^{-1}(\log n)}{c_n} \hbox{ i.o. }} =
\left\{%
\begin{array}{l}
    0 \\
    1
\end{array}%
\right. \Longleftrightarrow
\sum_{n\geq 2} \frac{1}{n\log n(c_n)^2}
\left\{%
\begin{array}{l}
    <\infty \\
    =\infty.
\end{array}%
\right.
\end{equation*}
In particular, with probability $1$:
\begin{equation*}
\liminf_{n\to\infty}\frac{\Pare{\log\log
n}^{\beta}}{a^{-1}\Pare{\log n}}\sup_{k\leq n}|Z_k| =
\left\{%
\begin{array}{ll}
    0, & \hbox{if $\beta\leq 1/2$}, \\
    \infty, & \hbox{if $\beta>1/2$}. \\
\end{array}%
\right.
\end{equation*}
\end{theop}

\textbf{Acknowledgements}: I would like to thank my Ph.D.
supervisor Yueyun Hu for his help and advices.
\bibliographystyle{plain}
\bibliography{biblio}
\end{document}